\documentclass[a4paper,11pt]{amsart}

\usepackage{array}
\usepackage{booktabs}
\usepackage{multicol}
\usepackage{multirow}
\usepackage{tabularx}
\usepackage{subfigure}
\usepackage{tikz}
\usepackage{cite}
\usepackage{soul}

\usepackage{systeme}
\usepackage{graphicx}
\usepackage[]{units}
\usepackage{color}
\usepackage{mathrsfs}
\usepackage{bm}
\usepackage[ansinew]{inputenc}
\usepackage[breaklinks,hidelinks]{hyperref}
\usepackage{float}
\usepackage{amsmath}
\usepackage{isomath}
\usepackage{amsthm}
\usepackage{amssymb}
\usepackage{amssymb}

\usepackage[foot]{amsaddr}
\usepackage[margin=2.7cm]{geometry}

\theoremstyle{remark}

\usepackage{hyperref}

\newcommand{\ldb}{\mathopen{\lbrack\!\lbrack}}
\newcommand{\rdb}{\mathclose{\rbrack\!\rbrack}}

\usepackage{bm}
\usepackage{graphicx,color, amsfonts,amsmath}

\usepackage{mwe}

\DeclareSymbolFont{matha}{OML}{txmi}{m}{it}
\DeclareMathSymbol{\varv}{\mathord}{matha}{118}

\begin{document}

\title[Geometrically parameterized Navier--Stokes CutFEM ROMs]{A Reduced Order Cut Finite Element method for geometrically parameterized steady and unsteady Navier--Stokes problems}
\author{Efthymios N. Karatzas\textsuperscript{1,4}}
\email{karmakis@math.ntua.gr}

\author{Monica Nonino\textsuperscript{2}}
\email{monica.nonino@univie.ac.at}

\author{Francesco Ballarin\textsuperscript{3}}
\email{francesco.ballarin@unicatt.it}

\author{Gianluigi Rozza\textsuperscript{4}}
\email{grozza@sissa.it}

\address{\textsuperscript{1}NTUA, School of Applied Mathematical and Physical Sciences, Department of Mathematics, Athens  and FORTH Institute of Applied and Computational Mathematics,
Heraclion, Crete, Greece.}
\address{\textsuperscript{2}University of Vienna, Department of Mathematics, Vienna, Austria.}
\address{\textsuperscript{3}Catholic University of the Sacred Heart, Department of Mathematics, Milan, Italy.}
\address{\textsuperscript{4}SISSA, International School for Advanced Studies, Mathematics Area, mathLab, Trieste, Italy.}

\subjclass[2010]{78M34, 97N40, 35Q35}

\begin{abstract}
We focus on steady and unsteady Navier--Stokes flow systems in a reduced-order modeling framework based on Proper Orthogonal Decomposition within a levelset geometry description and discretized by an unfitted mesh Finite Element Method. This work extends the approaches of  \cite{KaBaRO18,KaratzasStabileNouveauScovazziRozza2018,KaratzasStabileNouveauScovazziRozzaNS2019} to nonlinear CutFEM discretization. 
We construct and investigate a unified and geometry independent reduced basis which overcomes many barriers and complications of the past, that may occur whenever geometrical morphings are taking place. By employing a geometry independent reduced basis, we are able to avoid remeshing and transformation to reference configurations, and we are able to handle complex geometries. This combination of a fixed background mesh in a fixed extended background geometry with reduced order techniques appears beneficial and advantageous in many industrial and engineering applications, which could not be resolved efficiently in the past.

\end{abstract}
\maketitle
\section{Introduction}\label{sec:intro}
In the present work, we are interested in studying geometrically parameterized steady and unsteady Navier--Stokes equations in a Eulerian framework. We rely on an approach based on unfitted mesh Finite Element Method, 
which shows its flexibility especially when domains are subject to large deformations, and classical methods such as the Finite Element Method (FEM) fail.

In general, new computational tools have been invented and studied over the past years to solve numerically Navier-Stokes problems: the classical FEM is a powerful tool to discretize the physical domain of interest and simulate the behavior of the solution, and its efficiency has been proven in a wide range of applications \cite{boffi_mixed, Richter, NoBaRo19}. Nonetheless, FEM capability to handle geometrically parametrized problems comes to a limit, this limit being given by extremely complex geometries, but also by situations where large deformations, fractures, contact points occur. As an alternative to classical FEM, we can consider Finite Element (FE) approximations of the physical fields that are not fitted to the actual physical geometry. The FE approximations are then cut at the boundaries and interfaces: this gives rise to the Cut Finite Element Method (CutFEM). For a more precise idea and for more rigorous definitions of what ``cutting'' a physical field means, and for a detailed introduction to CutFEM, we refer to  \cite{burman_massing_hansbo} and references therein.\\
The repeated solution of parametrized problems discretized by CutFEM on the other hand, is an expensive task (whose cost essentially depends on the size of the underlying background mesh), especially in complex geometries.
It is precisely at this point that the Reduced Basis Method (RBM) comes into play.
It is well known that the Reduced Basis Method is an extremely powerful tool to obtain a speedup in the simulation of the behavior of the solution of the system. The method relies on a set of already computed solutions (snapshots) for different parameter values: see, for example, \cite{HeRoSta16,RozzaHuynhPatera1,HaasdonkOhlberger}. Therein, these snapshots are FE approximations of the truth solution, thus the RBM relies on the FEM.
{Even though in general there are several methods that can be employed to project the full order system to a reduced system, see for example \cite{Rozza2008229,ChinestaEnc2017,Kalashnikova_ROMcomprohtua,quarteroniRB2016,Chinesta2011,Dumon20111387}, in the present work we will employ the Proper Orthogonal Decomposition (POD).
We use a fixed background geometry and mesh: this approach leads to important advantages whenever a geometry deforms,\cite{KaBaRO18}, and it overcomes several related limitations in efficiency, compared with traditional FEM, see e.g.\cite{BaFa2014,karatzas_stabile_atallah,KaratzasStabileNouveauScovazziRozza2018}. For the reader interested in Reduced Order Methods based on classical FEM  we refer to: \cite{HeRoSta16} for Proper Orthogonal Decomposition, to \cite{ChinestaEnc2017,HeRoSta16,Kalashnikova_ROMcomprohtua,quarteroniRB2016,Rozza2008229} for greedy approaches and certified Reduced Basis Method, to \cite{ChinestaEnc2017,Chinesta2011,Dumon20111387} for Proper Generalized Decomposition (PGD), to \cite{Rozza2008229,grepl2005} for linear elliptic and parabolic systems and to \cite{Veroy2003,Grepl2007,RoQua07} for nonlinear problems.}
\\
{The aim of this work is to implement a Reduced Basis Method for the stationary and time--dependent Navier--Stokes equations, that relies on an unfitted CutFEM discretization.  Already existing results on RBM applied to Navier Stokes problems with standard FEM can be found in the classical literature, see for example \cite{doi:10.1002/9781119176817.ecm2051, RoQua07,Stabile2019}; for other results concerning RBM based on embedded FEM for different kind of problems, see for example \cite{karatzas_stabile_atallah,KaratzasStabileNouveauScovazziRozza2018}, as well as \cite{KaratzasStabileNouveauScovazziRozzaNS2019} for Navier-Stokes with the Shifted Boundary Method.
Nevertheless, there is still a fundamental lack in the literature for results that combine the RBM with an unfitted CutFEM discretization: indeed, a first important reference in this framework is represented by \cite{KaBaRO18}, where the authors focus on steady, linear model problems such as the Stokes problem and the Darcy problem.
The results presented in this manuscript represent the first step in the application of unfitted CutFEM based Reduced Basis Method to time--dependent, nonlinear fluid flows: to the best of our knowledge, indeed, existing results in the reduced order model framework include only linear problems, see \cite{KaBaRO18}.}

{
The paper is structured as follows: in Section \ref{sec:HF} we introduce some basic notions of unfitted mesh Finite Element discretization, as well as some notation that will be used throughout the work.
In Section \ref{sec:Steady_NS} we introduce the steady Navier--Stokes equations with incompressibility constraint in a fluid domain where the shape of some of the boundaries is described through a levelset function depending on a geometrical parameter. In Section \ref{discrete steady Stokes weak formulation} we state the weak formulation of the problem of interest, which is based on Nitsche's method with penalty term, and a Ghost Penalty stabilization. In Section \ref{sec:ROM} we present the Proper Orthogonal Decomposition, with a focus in Section \ref{lifting function} on the lifting of non--homogeneous Dirichlet boundary conditions that are imposed strongly, and with a focus in Section \ref{natural smooth extension} on the \emph{natural smooth extension}, a technique that is here used to obtain improved parameter--independent reduced basis functions. In Section \ref{sec:online} we present the online system to be solved, and finally in Section \ref{numerical results steady} we show some numerical results.
In Section \ref{Unsteady NS} we move then to the unsteady Navier--Stokes equations, formulated over the same physical domain as the one previously considered: the strong formulation is presented in Section \ref{strong formulation unsteady}. In Section \ref{Time discretization} we state the weak formulation after time discretization and after spatial discretization. In Section \ref{POD unsteady} we present the POD technique used in the case of time--dependent parametrized problems, and in Section \ref{numer_examples_unsteady} we present some numerical results. In Section \ref{immersed obstacle} we introduce another computational fluid dynamics test case, where an obstacle is now immersed in a fluid. In Section \ref{immersed obstacle weak formulation} we briefly recall the weak formulation, as well as the POD and the online system in Section \ref{immersed obstacle POD}; numerical results are provided in Section \ref{immersed obstacle numerical results}.
Conclusions and perspectives are provided in Section \ref{sec:conclusions}.
}
\section{Full order discretization by CutFEM:  an introduction to terminology and definitions} \label{sec:HF} 
The aim of this Section is to introduce some basic notions and definitions in the CutFEM framework: these definitions will be employed in the discrete formulation of the problems that we present hereafter.

{As mentioned in the Introduction, unfitted Finite Element discretizations are extremely useful for the numerical simulation of problems whose physical domain undergoes a significant change in its topology: contact points occur, overlapping domains, changes in the shape due to geometrical parametrization are just some of many examples. The fundamental idea at the core of unfitted discretization is the realization of a \emph{fixed background mesh}: this mesh, once generated, will not change, thus avoiding any additional expensive procedures as remeshing. The physical domain over which the problem of interest is formulated will intersect the elements of the background mesh, creating a natural subdivision into two meshes: an active mesh and an inactive mesh; in addition, the boundary of the physical domain will cut, in an arbitrary way, some of the elements of the background mesh. All these concepts are at the basis of every unfitted FEM discretization: let us now see more in detail how all these entities are defined.}\\
\begin{figure}
\centering
\begin{tikzpicture}
\node[anchor=south west,inner sep=0] (image) at (0,0) {\includegraphics[scale=1.5]{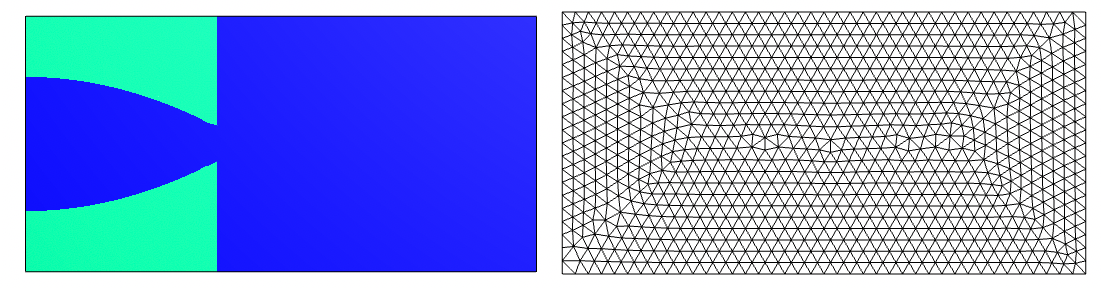}};
\begin{scope}[x={(image.south east)},y={(image.north west)}]
	\draw[white,thick] (0.25,0.5) node {{$\Omega$}};
	\draw[white,thick] (0.05,0.5) node {{$\Gamma_{in}$}};
	\draw[white,thick] (0.35,0.85) node {{$\Gamma_D$}};
	\draw[white,thick] (0.35,0.15) node {{$\Gamma_D$}};
	\draw[white,thick] (0.45,0.5) node {{$\Gamma_{out}$}};
	\draw[black,thick] (0.1,0.75) node {{$\Gamma$}};
	\draw[black,thick] (0.1,0.25) node {{$\Gamma$}};
	\draw[black,thick] (0.7,0.5) node {{$\hat{\mathcal{I}}_h$}};
\end{scope}
\end{tikzpicture}
\caption{Left: the physical domain of interest $\Omega$ (dark blue), with inflow boundary $\Gamma_{in}$, outflow boundary $\Gamma_{out}$, and top and bottom boundaries which are denoted here $\Gamma_D$. The boundary $\Gamma$ is the boundary delimited by the green sets, and is the immersed boundary. Right: the background mesh $\hat{\mathcal{I}}_h$.} \label{background mesh}
\end{figure}
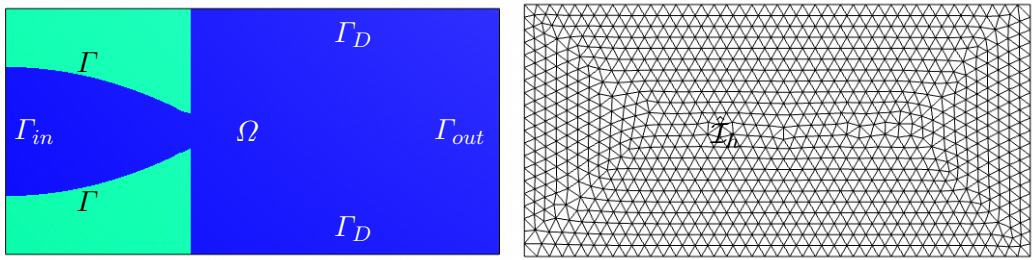

\begin{figure}
\centering
\begin{tikzpicture}
\node[anchor=south west,inner sep=0] (image) at (0,0) {\includegraphics[scale=1.5]{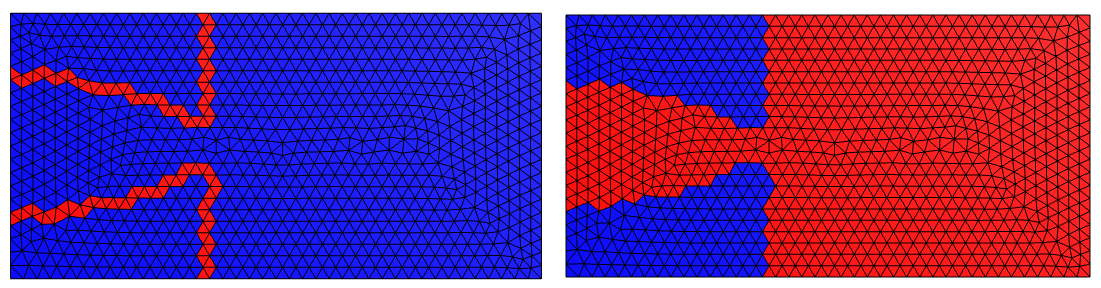}};
\begin{scope}[x={(image.south east)},y={(image.north west)}]
	\draw[white,thick] (0.25,0.5) node {{$\mathcal{I}_{\Gamma}$}};
	\draw[white,thick] (0.8,0.5) node {{${\mathcal{I}}_h$}};
\end{scope}
\end{tikzpicture}
\caption{Left: the portion of the mesh intersected by the immersed boundary $\Gamma$, namely $\mathcal{I}_{\Gamma}$, is depicted in red. Right: the fictitious domain $\Omega^{*}_h$ is depicted in red, and is determined by the elements $T$ of the active mesh $\mathcal{I}_h$.}\label{cut elements}
\end{figure}

Let $\Omega$ be the physical domain over which our problem is formulated, and let $\hat{\mathcal{I}}_h$ be a fixed background mesh of our choice, of mesh size $h>0$, covering $\Omega$, {{ with inlet boundary $\Gamma_{in}=\{\bm{x}\in\Gamma:\bm{u}(\bm{x},t)\cdot\bm{n}<0\}$}} (see Figure \ref{background mesh}). Given a background mesh, this identifies in a natural way the \emph{active mesh} $\mathcal{I}_h$, which is the portion of the background mesh made by elements $T$ that are actually intersected by the physical domain:
\begin{equation*}
\mathcal{I}_h := \{T \in \hat{\mathcal{I}}_h \colon \Omega \cap T \neq \emptyset\}.
\end{equation*}
Then, starting from the active mesh we can define a domain $\Omega_h^{*}$ as follows:
\begin{equation*}
\Omega_h^{*} = \bigcup\limits_{T \in \mathcal{I}_h}^{} T.
\end{equation*}
Two situations can present: $\mathcal{I}_h$ is a
\begin{itemize}
\item \emph{boundary fitted mesh} if $\overline{\Omega^{*}_h} = \overline{\Omega}
$;
\item \emph{unfitted mesh} if $\overline{\Omega}
\subsetneq \overline{\Omega^{*}_h}$. In this case $\Omega^{*}_h$ is the so called \emph{fictitious domain}, and we have an unfitted discretization of the problem of interest (see Figure \ref{cut elements}). This is the kind of discretization that will adopted throughout the paper.
\end{itemize}
For unfitted meshes a crucial role is played by those elements $T$ of the active mesh that are cut by the immersed boundary $\Gamma$ of the physical domain $\Omega$:
\begin{equation*}
\mathcal{I}_{\Gamma} := \{T\in \mathcal{I}_h \colon T\cap \Gamma\neq\emptyset\}.
\end{equation*}
Related to the set $\mathcal{I}_{\Gamma}$ we can also define the set of facets that belong to elements intersected by the boundary:
\begin{equation*}
\mathcal{F}_{\Gamma} :=\{F \in \mathcal{F}\colon F \text{ is a facet of an element } T\in \mathcal{I}_{\Gamma}\}.
\end{equation*}

{In the framework of unfitted meshes and unfitted Finite Element discretization, the role of the cut elements $T\in \mathcal{I}_{\Gamma}$ is important: indeed, it has been shown (see for example \cite{burman_massing_hansbo, burman_hansbo, BuHa11, BURMAN20101217}) that stability issues may arise, depending on the quality of the cut, namely depending on how the boundary $\Gamma$ cuts an element $T$. In order to overcome the dependency of the stability and a priori estimates on the position of the interface, and the overall ill-conditioning of the global system matrix due to bad intersections, Burman et al. introduced a stabilization technique called Ghost Penalty, see \cite{burman_massing_hansbo,burman_hansbo,burman_fernandez,BURMAN20074045}. The idea of the Ghost Penalty stabilization is to adding weakly consistent operators with the aim of having a better control on the solution in $\overline{\Omega_h^{*}}\setminus\overline{\Omega}$. We will return more into the details of the Ghost Penalty stabilization used in Section \ref{discrete steady Stokes weak formulation}. We remark that the Ghost Penalty technique is the stabilization that will be used throughout the paper, nevertheless there are alternative techniques like the residual-based stabilization (RBVM), see for example \cite{schott2017stabilized}}. 

\section{Steady Navier--Stokes}\label{sec:Steady_NS}
In the following, we focus on the steady incompressible Navier--Stokes equations, which are formulated within an Eulerian formalism. We first introduce the problem formulation, and then we present the discretized version of the original problem, in the CutFEM framework.
\begin{figure}
\centering
\begin{tikzpicture}
\node[anchor=south west,inner sep=0] (image) at (0,0) {\includegraphics[scale=0.7]{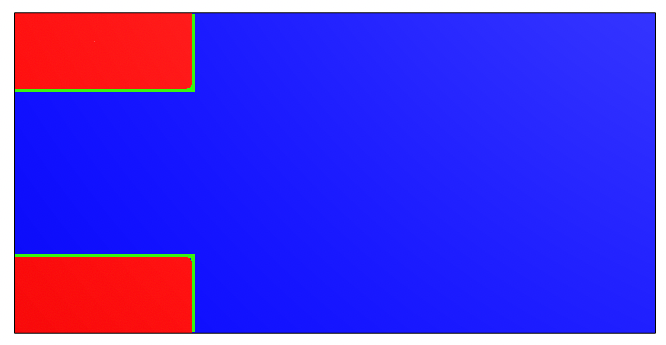}};
\begin{scope}[x={(image.south east)},y={(image.north west)}]
	\draw[white,thick] (0.08,0.5) node {$\Gamma_{in}$};
	\draw[white,thick] (0.88,0.5) node {$\Gamma_{out}$};
	\draw[white,thick] (0.5,0.9) node {$\Gamma_{D}$};
	\draw[white,thick] (0.5,0.1) node {$\Gamma_{D}$};
	\draw[white,thick] (0.3,0.7) node {$\Gamma_{\theta}$};
	\draw[white,thick] (0.2,0.85) node {$\mathcal{D}(\theta)$};
	\draw[white,thick] (0.2,0.15) node {$\mathcal{D}(\theta)$};
\end{scope}
\end{tikzpicture}
\caption{Instance geometry of the problem for $\theta=0.37$: the fluid domain $\Omega(\theta)$ (blue) has an inlet boundary $\Gamma_{in}$ on the left and an outlet boundary $\Gamma_{out}$ on the right. The rest of the boundaries ($\Gamma_D$) are Dirichlet type boundaries.
The shape of $\mathcal{D}(\theta)$ (red) is described by the levelset function $\Phi_{\theta}(x,y)$ (see Equation \ref{levelset_equation}). The immersed boundary $\Gamma_{\theta}$ is depicted in green.
}
\label{mesh}
\end{figure}
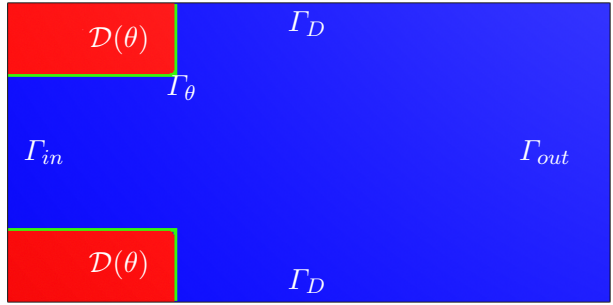
\subsection{Strong formulation} 
Let $\mathcal{R}$ be a background rectangular domain in $\mathbb{R}^2$, and let $\mathcal{D}(\theta)\subset\mathcal{R}$ be a bounded subset of $\mathcal{R}$, whose boundary is described through a levelset function $\{\Phi_{\theta}=0\}$, where $\Phi_{\theta}$ is an implicit function depending on a geometrical parameter $\theta$. The physical domain over which our problem is formulated is $\Omega(\theta) := \mathcal{R}\setminus\mathcal{D}(\theta)$, and is depicted in Figure \ref{mesh} for a given value of $\theta$. {In Figure \ref{mesh} we can see the inlet boundary $\Gamma_{in}$,  the outlet boundary $\Gamma_{out}$, as well as the top and bottom wall of the domain, which are denoted by $\Gamma_D$, since Dirichlet conditions will be applied there. We further denote $\Gamma_{\theta}$ the immersed boundary, that is, the remaining part of the boundary $\partial\Omega(\theta)$ that is in common to $\mathcal{D}(\theta)$.}
We denote by $\mathcal{P}$ the parameter space to which $\theta$ belongs.
Under these assumptions, our problem of interest reads: for every $\theta\in\mathcal{P}$, find $\bm{u}(\theta)\colon\Omega(\theta)\mapsto\mathbb{R}^2$ and $p(\theta)\colon\Omega(\theta)\mapsto\mathbb{R}$
such that:
\begin{equation}\label{eq:NS_strong}
\begin{cases}
-\mu\Delta \bm{u}(\theta) + \nabla p(\theta) + (\bm{u}(\theta)\cdot \nabla)\bm{u}(\theta) =  \bm{f}(\theta) \quad \text{in $\Omega(\theta)$}, \\
\text{div}\bm{u}(\theta)=0 \quad \text{in $\Omega(\theta)$},
\end{cases}
\end{equation}
where 
$\mu$ is the fluid viscosity and $\bm{f}(\theta)$ is the fluid volume external force. {Problem \eqref{eq:NS_strong} is completed by the following boundary conditions: at the inlet boundary $\Gamma_{in}$ we impose a prescribed inlet velocity (througout the manuscript we choose a parameter independent inlet profile) $\bm{u}(\theta) = \bm{u}_{in}$; we then have a zero outflow condition on $\Gamma_{out}$, a no--slip boundary condition $\bm{u}(\theta)\cdot\bm{n}_D=0$ on $\Gamma_D$, with $\bm{n}_D$ the outgoing normal to $\Gamma_D$, and a no--slip boundary condition $\bm{u}(\theta) = \bm{0}$ on $\Gamma_{\theta}$.}
\subsubsection{Geometrical parametrization}
For the problem considered in this section the expression of the levelset function is the following:

\begin{equation}
\begin{aligned}
\Phi_{\theta}(x, y)= {} &
-\Bigl(\left|A(x) + B(x)-1\right|
 +\left|{A(x)- B(x)-2}\right| + D(x)\Bigr) \\
& \cdot\Bigl(\left|{A(x) + C(x)-1}\right| +\left|{A(x)-C(x)-2}\right| + D(x)\Bigr),
\end{aligned}
\label{levelset_equation}
\end{equation}
where $A(x)=\sqrt{k_1}\left|x-k_3\right|$, $B(x)=\sqrt{k_2}\left|y-k_4\right|$, $C(x)=\sqrt{k_2}\left|y-k_5\right|$ and $D(x)=e^{-\theta}(k_1(x-k_3)^2)\theta-4$.
\begin{figure}
\centering
\includegraphics[scale=0.221]{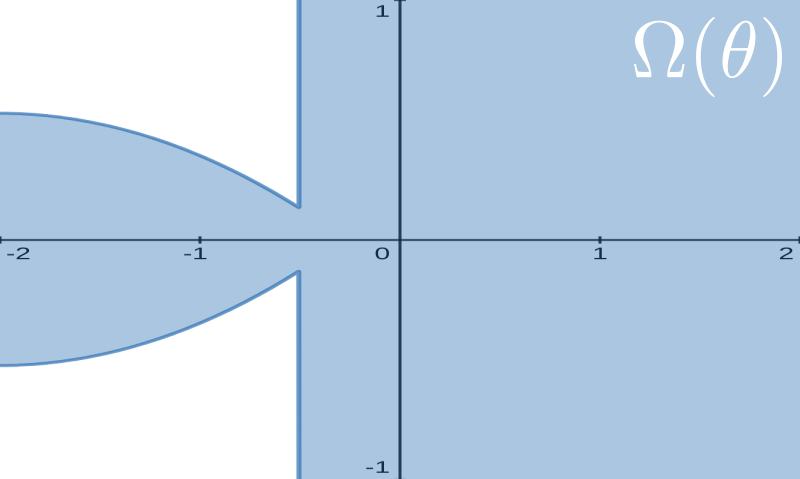}\qquad
\includegraphics[scale=0.221]{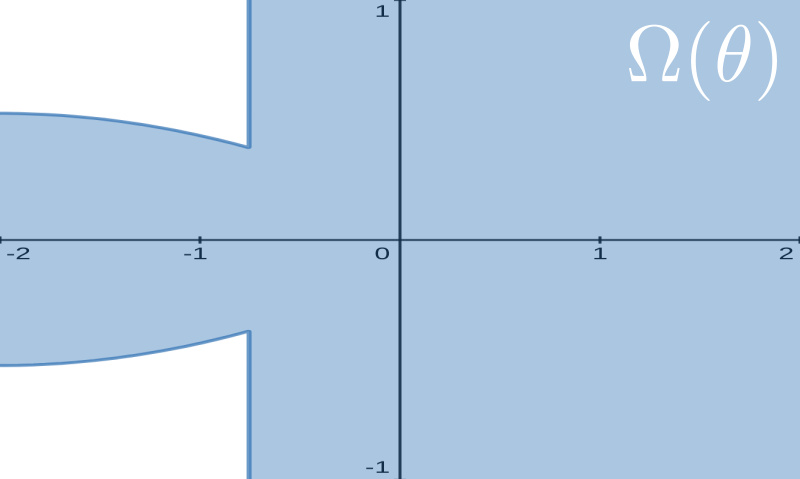}\qquad
\includegraphics[scale=0.221]{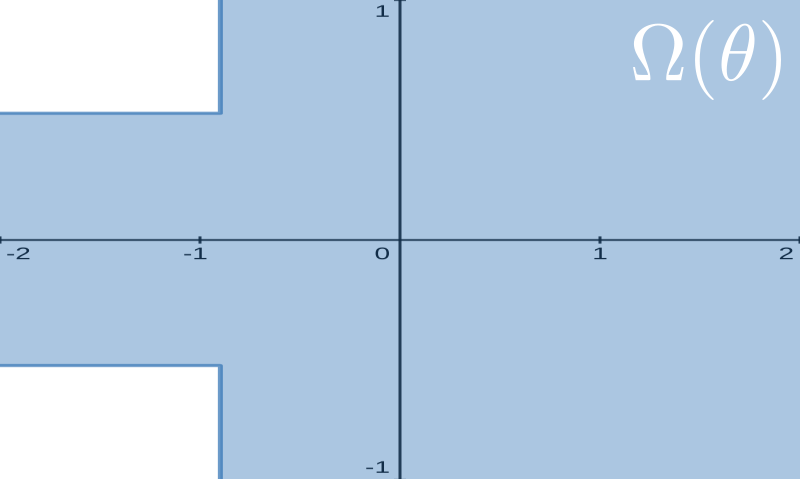}\\
\vskip10pt
\includegraphics[scale=0.221]{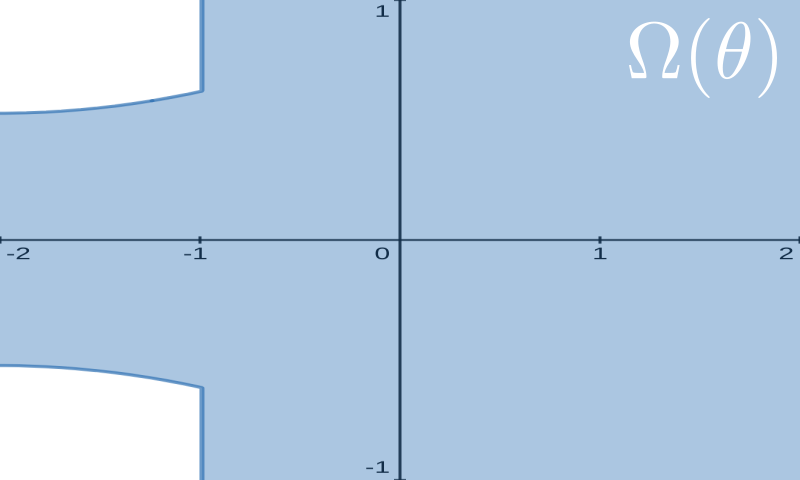}\qquad
\includegraphics[scale=0.221]{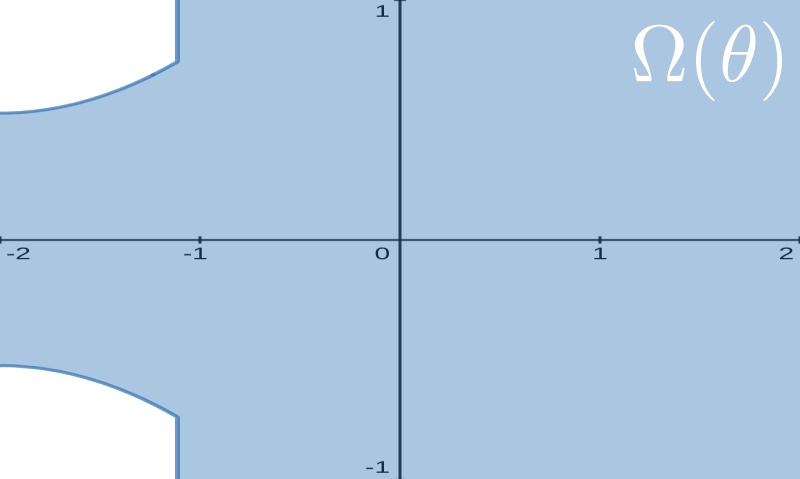}\qquad
\includegraphics[scale=0.221]{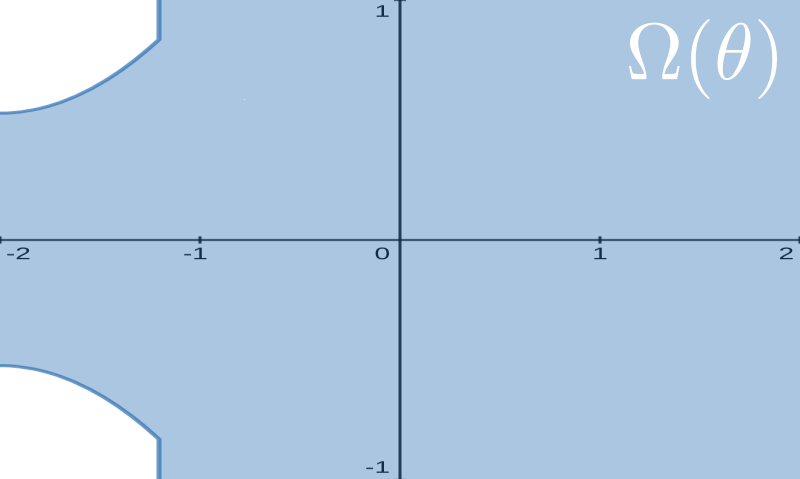}
\caption{Six examples of solid walls, described by the levelset $\{\Phi_{\theta}=0\}$. From left to the right, the levelset for $\theta=-0.1, -0.06, 0, 0.18, 0.37, 0.50$.}
\label{levelset_image}
\end{figure}

The values of the constants $k_1, k_2,k_3, k_4, k_5$ are reported in Table \ref{levelset parameter values}.
To have a better idea of how the shape of the walls changes by varying the parameter $\theta$, the reader is referred to Figure \ref{levelset_image}.
\begin{table}
\centering
\caption{Values for the constants in the levelset equation \eqref{levelset_equation}.}
\label{levelset parameter values}
\begin{tabular}{|l|c|}
\hline
\textbf{Constant} & \textbf{Value}\\
\hline
$k_1$ & $10$\\
$k_2$ & $10$\\
$k_3$ & $-2$\\
$k_4$ & $-1$\\
$k_5$ & $1$\\
\hline
\end{tabular}
\end{table}

\subsection{Discrete weak formulation}\label{discrete steady Stokes weak formulation}
As we can see from Figure \ref{mesh}, the background mesh $\hat{\mathcal{I}}_h$ is a rectangular mesh made by triangular elements. By choosing an unfitted mesh, once we have defined a background mesh, there is no need to remesh every time the parameter $\theta$ changes (and hence every time the shape of the levelset in dark grey in Figure \ref{mesh} changes).
{Before going any further, we remark that we decide to impose the Dirichlet boundary conditions in the following way: the no--slip condition on the immersed boundary $\Gamma_{\theta}$ is imposed weakly, whereas the inflow condition, as well as the homogeneous Dirichlet condition on $\Gamma_D$ are imposed strongly.}
Let us now introduce the following discrete approximation spaces:
\begin{equation*}
\begin{split}
V_{h,{{2}}}(\theta) &:= \{\bm{v}_h \in C_0(\Omega^*_h(\theta))^2 \colon \bm{v}_h |_{T} \in (\mathcal{P}^{{2}}(T))^2, \quad\forall T \in \mathcal{I}_h(\theta)\},\\
Q_{h,1}(\theta) &:= \{q_h \in C_0(\Omega^*_h(\theta))\colon
q_h|_{T} \in \mathcal{P}^1(T), \quad\forall T \in \mathcal{I}_h(\theta)\},
\end{split}
\end{equation*}
where $\mathcal{I}_h(\theta)$ and $\Omega_h^*(\theta)$ are respectively the active mesh and the fictitious domain corresponding to the physical domain $\Omega(\theta)$, as defined in Section $2.1$: these spaces clearly depend on $\theta$, since the levelset geometry changes according to the parameter.
{We then define:
\begin{equation*}
\begin{split}
V_{h,{{2}}}^D\theta) &:= \{\bm{v}_h \in V_{h,{{2}}}(\theta) \colon \bm{v}_h=\bm{u}_{in}\text{ on }\Gamma_{in}\text{ and }\bm{v}_h\cdot\bm{n}_D=0\text{ on }\Gamma_D\},\\
V_{h,{{2}}}^0(\theta) &:= \{\bm{v}_h \in V_{h,{{2}}}(\theta) \colon \bm{v}_h=\bm{0}\text{ on }\Gamma_{in}\}
\end{split}
\end{equation*}
}
The discretized weak formulation of the problem, with Nitsche's Method and Ghost Penalty terms then reads: find $(\bm{u}_h(\theta), p_h(\theta)) \in V_{h,2}^D(\theta)\times Q_{h,1}(\theta)$
such that for all test functions $(\bm{v}_h(\theta), q_h(\theta))\in V_{h,2}^0(\theta)\times Q_{h,1}(\theta)$:
\begin{equation}
\label{weak formulation}
\mathcal{A}({\bm{u_h}}(\theta),p_h({\theta}); \bm{v}_h(\theta),q_h(\theta)) = \mathcal{L}(\bm{v}_h(\theta)),
\end{equation}
where:
{
\begin{eqnarray*}
\mathcal{A}(\bm{u}_h(\theta),p_h(\theta); \bm{v}_h(\theta), q_h(\theta)) &=& {a_h(\bm{u}_h(\theta), \bm{v}_h(\theta)) + b_h(p_h(\theta), \bm{v}_h(\theta)) - b_h(q_h(\theta), \bm{u}_h(\theta)) +}\nonumber\\
&&{+ c_h(\bm{u}_h(\theta); \bm{u}_h(\theta), \bm{v}_h(\theta) + g^{GP}(\bm{u}_h(\theta), \bm{v}_h(\theta); p_h(\theta), q_h(\theta))}\nonumber\\
\mathcal{L}(\textbf{v}_h) &=& 
(f,v_h)_{\Omega(\theta)}
\end{eqnarray*}
}
In the previous equation we have the following bilinear (and trilinear) forms:
{
\begin{eqnarray*}
a_h(\bm{u}_h(\theta), \bm{v}_h(\theta))&=&{ \mu(\nabla\bm{u}_h(\theta), \nabla\bm{v}_h(\theta))_{\Omega(\theta)} -\mu(\nabla\bm{u}_h(\theta)\bm{n}_{\Gamma}, \bm{v}_h(\theta))_{\Gamma_{\theta}} -}\nonumber\\
&&{-\mu(\bm{u}_h(\theta), \nabla\bm{v}_h(\theta)\bm{n}_{\Gamma})_{\Gamma_{\theta}} +{\gamma\mu/}{h}(\bm{u}_h(\theta), \bm{v}_h(\theta))_{\Gamma_{\theta}}}
\nonumber\\
&&
{{+ {\gamma\phi}/{h}(\bm{u_h}(\theta)\cdot \bm{n}_{\Gamma},\bm{v_h}(\theta)\cdot \bm{n}_{\Gamma})_{\Gamma_{\theta}}}}\nonumber\\
c_h(\bm{u}_h(\theta); \bm{u}_h(\theta), \bm{v}_h(\theta))
&=&{((\bm{u}_h(\theta)\cdot\nabla)\bm{u}_h(\theta),\bm{v}_h(\theta))_{\Omega(\theta)}} -{((\bm{u}_h(\theta)\cdot\bm{n}
)\bm{u}_{in},\bm{v}_h(\theta))_{\Gamma_{in}}}
\nonumber\\
b_h(p_h(\theta), \bm{v}_h(\theta))&=& -{(p_h(\theta), \nabla\cdot\bm{v}_h(\theta))_{\Omega(\theta)} + (p_h(\theta)\bm{n}_{\Gamma}, \bm{v}_h(\theta))_{\Gamma_{\theta}}}\nonumber\\
g^{GP}(\bm{u}_h(\theta), \bm{v}_h(\theta); p_h(\theta), q_h(\theta))&=&
{g_u^{GP}(\bm{u}_h(\theta);\bm{u}_h(\theta), \bm{v}_h(\theta)) + g^{GP}_p(\bm{u}_h(\theta); p_h(\theta), q_h(\theta))}\nonumber\\
&&
{{+g_{\mu}^{GP}(\bm{u}_h(\theta), \bm{v}_h(\theta)) + g_{\beta}^{GP}(\bm{u}_h(\theta);\bm{u}_h(\theta), \bm{v}_h(\theta))}}\nonumber
\end{eqnarray*}s
where $\bm{n}_{\Gamma}$ is the outward pointing normal to the immersed boundary $\Gamma_{\theta}$, and $h$ is the mesh size.
The term $-((\bm{u}_h(\theta)\cdot\bm{n}
)\bm{u}_h(\theta), \bm{v}_h(\theta))_{\Gamma_{in}
} $ in the trilinear form $c_h$, where $\bm{n}
$ in this term is the outgoing normal to the inlet boundary, that accounts for the inlet flow $\Gamma_{in}=\{\bm{x}\in\Gamma:\bm{u}(\bm{x},t)\cdot\bm{n}<0\}$ (see for example \cite{https://doi.org/10.1002/cnm.529}). 
The ghost penalties 
$g_u^{GP}(\bm{u}_h(\theta), \bm{v}_h(\theta))$,  $g^{GP}_p(p_h(\theta), q_h(\theta))$, $g_\mu^{GP}(\bm{u}_h(\theta), \bm{v}_h(\theta))$,  and $g^{GP}_\beta(p_h(\theta), q_h(\theta))$ are defined as:
\begin{eqnarray}
g_u^{GP}(\bm{u}_h(\theta);\bm{u}_h(\theta), \bm{v}_h(\theta))&=&
{{{\gamma_u\sum_{F\in\mathcal{F}_{{\Gamma}_{\theta}}}\phi_{u, F,\theta} h^{2j+1}(\ldb\nabla\cdot\partial_n^j\bm{u}_h(\theta)\rdb, \ldb\nabla\cdot\partial_n^j\bm{v}_h(\theta)\rdb)_F}}},
\nonumber\\
g^{GP}_p(\bm{u}_h(\theta); p_h(\theta), q_h(\theta))&=& 
\gamma_p \sum_{F\in\mathcal{F}_{\Gamma_\theta}
}{h^3}{\mu}^{-1}({1}/{\text{max}({h}/{\mu}\lvert\lvert {\bm u}(\theta)\rvert\rvert_{\infty,T}, 1))} (\ldb\partial_n p_h(\theta)\rdb,\ldb\partial_n q_h(\theta)\rdb)_F,
\nonumber
\nonumber\\
{{g_{\mu}^{GP}(\bm{u}_h(\theta), \bm{v}_h(\theta))}} &=& 
{{{\gamma_{\mu}\sum_{j=1}^2\sum_{F\in\mathcal{F}_{{\Gamma}_{\theta}}}\mu h^{2j-1}(\ldb\partial_n^j\bm{u}_h(\theta)\rdb, \ldb\partial_n^j\bm{v}_h(\theta)\rdb)_F}\nonumber}}, \\
{{g_{\beta}^{GP}(\bm{u}_h(\theta);\bm{u}_h(\theta), \bm{v}_h(\theta))}}&=&{{{\gamma_{\beta}\sum_{j=1}^2\sum_{F\in\mathcal{F}_{{\Gamma}_{\theta}}}\phi_{\beta, F,\theta} u_{\infty, F,\theta}^2h^{2j-1}(\ldb(\bm{u}_h(\theta)\cdot\nabla)\partial_n^j\bm{u}_h(\theta)\rdb, \ldb(\bm{u}_h(\theta)\cdot\nabla)\partial_n^j\bm{v}_h(\theta)\rdb)_F}}}.\nonumber
\end{eqnarray}
In the previous equations, $
\mu
$ is the fluid kinematic viscosity, $\partial_n$ denotes the partial derivative with respect to the normal outgoing the face $F$. We have used the following notation:
\begin{eqnarray}
&& u_{\infty, F,\theta}:={\lvert\lvert \bm{u}_h(\theta)\rvert\rvert_{0, \infty, F}}, \quad
\label{u_infty}
\phi_{T,\theta}:={\mu + c_u\lvert\lvert\bm{u}_h(\theta)\rvert\rvert_{0, \infty, T}h_T}
\\
&&\phi_{\beta, T,\theta}:={\phi_{p, T,\theta}=h_T^2\phi_{T,\theta}^{-1}}, \quad
\phi_{u, T,\theta}:={\phi_{T,\theta}},\label{ultimo}
\end{eqnarray}
where $h_T$ is the characteristic element length.
Then, $\phi_{\beta, F,\theta}$, $\phi_{u, F,\theta}$ and $\phi_{p, F,\theta}$ are the corresponding face averages. In order to simplify the implementation and the exposition, we choose here $c_u=1.0$.
}
{The terms appearing in the bilinear form $a_h(\bm{u}_h, \bm{v}_h)$, in the bilinear form $b_h(p_h(\theta), \bm{v}_h(\theta))$ and in the trilinear form $c_h(\bm{u}_h(\theta; \bm{u}_h(\theta), \bm{v}_h(\theta))$ come from standard integration by parts of the steady Navier--Stokes equations, due to the fact that the test functions are non--vanishing on the boundary $\Gamma_{\theta}$. We then have the Nitsche terms, that have been used to weakly impose the non--slip boundary condition on $\Gamma_{\theta}$.
The term $g_{\beta}^{GP}$ is a convective stabilization term which for low Reynolds number is taken $\gamma_\beta=0$, whereas 
$g_{u}^{GP}$ is the incompressibility stabilization while $g_{p}^{GP}$, $g_{\mu}^{GP}$ contain additional terms that extend the solution to the extended domain, \cite{schott2017stabilized, burman_hansbo, KaBaRO18}. 
We point out that we choose a symmetric Nitsche method because, even if the non--symmetric alternative leads to better stability behaviour, it can also lead to suboptimal convergence or larger $L^2$ errors: for more details we refer to the discussions in \cite{penalty-free-burman, SchWa14}. For sake of completeness we report that the choice of a symmetric Nitsche method leads to some lack of coercivity in the bilinear form $a_h$, and therefore a stabilization term is consistently added. 
Additionally, the need for the Ghost Penalty stabilization terms is given by the fact that we are implementing an unfitted mesh discretization of the problem. With unfitted meshes indeed, it may happen that the cut elements (namely the elements $T\in\mathcal{I}_{\Gamma_{\theta}}$) are cut in an arbitrary way by the boundary $\Gamma_{\theta}$. The abritrariness of the cut position can lead to stability issues that ghost penalty stabilization terms can handle. In this manuscript we adopted the stabilization presented, for example, in \cite{schott2017stabilized}, as well as for an analysis on the stability properties of the chosen weak formulation.
 The values of the constants appearing in the previous equations as well as for the supremizer enrichment are reported in Table \ref{parameter values}.
 \begin{table}
\centering
\caption{Constants values for the weak formulation of the problem}
\label{parameter values}
\begin{tabular}{|l|c|l|c|}
\hline
\textbf{Constant} & \textbf{Value} & \textbf{Constant} & \textbf{Value}\\
\hline
$\mu$ & $0.05$ & $h$ & $0.07$\\
$\gamma$ & $10$ & $\gamma_u$ & $0.001$\\
$\alpha$ & $0.1$ & $\gamma_\mu$ & $0.1$\\
$\gamma_p$ & $0.1$ & $\gamma^0_s$ & $0.1$\\
$\lambda_s$ & $10$ & $\gamma^1_s$ & $0.01$\\
\hline
\end{tabular}
\end{table}
}
\subsection{Proper Orthogonal Decomposition-Galerkin Model Reduction}\label{sec:ROM}
We now present the Proper Orthogonal Decomposition (POD), a method that is applied to parameter-dependent problems in order to extract and create a set of reduced basis functions, with which we will subsequently perform the model order reduction. The POD method consists of two phases: one \textit{offline} and one \textit{online}. During the offline phase, we compute the solution of the problem of interest, for different values of the parameter $\theta$. These values of the parameter are collected from a training set $\mathcal{P}_{train}$, and the corresponding solutions are stored into a matrix, the so-called snapshots matrix. This matrix is then processed in order to extract the reduced basis functions.
Afterwards, in the online phase, we employ these basis functions in a way that reduces the dimension of the original problem, and in a way that is computationally efficient for (in our case) geometrically parametrized systems.
We remind hereafter that for POD-Galerkin ROMs for incompressible Navier--Stokes equations, instabilities in the approximation of the pressure may occur. We refer to \cite{Caiazzo2014598,Gerner2012,RoHuMa13} for a more detailed analysis of the problem, while for such instabilities on transient problems we refer to \cite{Iollo2000,Akhtar2009,Bergmann2009516,Sirisup2005218,taddei2017}.
For SUPG and PSPG kind of stabilization we refer to \cite{ballarin2015supremizer,RoVe07, stabile_stabilized, StaHiMoLoRo17}. Levelset techniques with cut Finite Element Nitsche method will be employed for the parametrization and ROM focusing on a fixed, geometrical parameter independent, background mesh following approaches as in \cite{KaBaRO18, KaratzasStabileNouveauScovazziRozza2018, karatzas_stabile_atallah}.
{
\subsubsection{The lifting function} \label{lifting function}
In view of the online phase of the reduction procedure, we have to take care of the non--homogeneous Dirichlet boundary condition at the inlet boundary $\Gamma_{in}$. We will now briefly explain how we overcome this difficulty, with the understanding that the following procedure will be carried out for all the test cases considered, even if we will not explicitly mention it again, for the sake of the simplicity of the notation.\\
The idea is that, in order to take care of a non--homogeneous Dirichlet boundary condition, we can define, for each parameter $\theta\in\mathcal{P}$, a lifting function $\bm{\ell}(\theta)\colon\Omega(\theta)\mapsto\mathbb{R}^2$, sucht that $\bm{\ell}(\theta)=\bm{u}_{in}$ on the inlet boundary $\Gamma_{in}$ and $\bm{\ell}(\theta)\cdot\bm{n}_D=0$ on $\Gamma_D$. Once we have done this, we can define a homogenized velocity 
\begin{equation*}
\bm{u}^0(\theta):=\bm{u}(\theta)-\bm{\ell}(\theta),
\end{equation*}
with $\bm{u}^0(\theta)\in {V}_{h, 2}^0(\theta)$.
This will become important in the next paragraph, when we state the algebraic formulation of the system that we are going to solve during the online phase of the method. The reader interested in more details on the lifting function and its use is referred, for example, to \cite{IAPICHINO2016408, ballarin2015supremizer}. }
{
\subsubsection{The natural smooth extension}\label{natural smooth extension}
Let us now make an important remark about the discrete function spaces and about the CutFEM solutions: as we can see from the definitions in Section \ref{discrete steady Stokes weak formulation}, the FE velocity $\bm{u}_h(\theta)$ belongs to a space that is $\theta$--dependent, namely $V_{h,2}(\theta)$, and the same goes for the FE pressure $p_h(\theta)\in Q_{h,1}(\theta)$. This $\theta$--dependence can create numerous problems: in general we can expect that, for two different values $\theta^1$ and $\theta^2$ of the parameter,  $\Omega(\theta^1)\neq \Omega(\theta^2)$, and therefore $\bm{u}_h(\theta^1)\notin V_{h,2}(\theta^2)$, and similarly for the pressure. This can lead to many difficulties, especially in view of the reduction procedure: in order to perform a Proper Orthogonal Decomposition we will have to compute scalar products between functions that, in theory, can belong to different discrete spaces.
In order to overcome this problem, the main idea is to extend the CutFEM solutions, namely the snapshots, to the whole background mesh: this is achieved thanks to a \emph{natural smooth extension} of both velocity and pressure. The realization, at the discrete level, of the natural smooth extension is carried out through the stabilization Ghost Penalty terms that appear in the weak formulation of the original problem: we refer the reader interested in more details and in the stability analysis of this procedure to \cite{BeBuHa09, refId0}. As we know, the CutFEM solution is defined up to $\Omega_h^*(\theta)$: as an effect of the ghost penalty stabilization, the solution smoothly goes to zero in the cut zone $\mathcal{I}_{\Gamma_{\theta}}$. Thanks to the ghost penalty stabilization then, we can simply extend the CutFEM solution to be zero in $\mathcal{R}\setminus\Omega_h^*(\theta)$, thus obtaining a solution that is defined on the whole background mesh.
Another possibility, which we do not take into consideration in this manuscript, is to implement instead an harmonic extension of the CutFEM solution: we refer the reader interested in a detailed discussion on techniques to extend the snapshots to the background mesh to \cite{KaBaRO18, BallarinRozzaMaday}: the reason we do not pursue this idea here is that the implementation of an harmonic extension would require an additional problem to be solved for every snapshot computed, thus incrementing the cost of the reduction procedure.\\
We remark that the natural smooth extension is employed also in order to extend the lifting function $\ell(\theta)$ to a function that is defined on the whole background mesh. Thanks to the extension procedure we obtain snapshots $(\hat{\bm{u}}^0_h(\theta), \hat{p}_h(\theta))$ that are defined on the common background mesh $\hat{\mathcal{I}_h}$, where now $\hat{\bm{u}}^0_h(\theta):=\hat{\bm{u}}_h(\theta)-\hat{\ell}(\theta)$.
Such extension defines a pair of velocity--pressure snapshots $(\hat{\bm{u}}^0_h(\theta), \hat{p}_h(\theta))$ belonging to \emph{$\theta$--independent} discrete spaces:
\begin{equation*}
\begin{split}
\hat{V}_{h,2}^0 &:= \{\hat{\bm{v}}_h \in C_0(\mathcal{R})^2 \colon \hat{\bm{v}}_h |_{T} \in (\mathcal{P}^2(T))^2 \text{ }\forall T \in \hat{\mathcal{I}}_h \text{ and }\bm{v}_h=\bm{0}\text{ on }\Gamma_{in} \text{ and }\bm{v}_h\cdot\bm{n}_D=0\text{ on }\Gamma_D\},\\
\hat{Q}_{h,1} &:= \{\hat{q}_h \in C_0(\mathcal{R})
\colon
\hat{q}_h |_{T} \in \mathcal{P}^1(T) \text{ }\forall T \in \hat{\mathcal{I}}_h\}.
\end{split}
\end{equation*}
At this point we can define a bijection between $\hat{V}_{h,2}^0$ and $\mathbb{R}^{\hat{N}_u^h}$ (respectively $\hat{Q}_{h,1}$ and $\mathbb{R}^{\hat{N}_p^h}$), where $\hat{N}_u^h$ and $\hat{N}_p^h$ are the dimensions of the discrete spaces $\hat{V}_{h,2}^0$ and $\hat{{Q}}_{h,1}$:
\begin{equation}
\label{bijection}
\begin{cases}
\underline{\textbf{u}}^0_h(\theta) = ({u}_h^1(\theta), \dots, u_h^{\hat{N}_u^h}(\theta))^T\in\mathbb{R}^{\hat{N}_u^h}\iff\hat{\bm{u}}^0_h(\theta)=\sum_{i=1}^{\hat{N}_u^h}u_h^i(\theta)\hat{\bm{\varphi}}^i\in \hat{V}_{h,2}^0,\\
\underline{p}_h(\theta) = ({p}_h^1(\theta), \dots, p_h^{\hat{N}_p^h}(\theta))^T\in\mathbb{R}^{\hat{N}_p^h}\iff p_h(\theta)=\sum_{i=1}^{\hat{N}_p^h}{p}_h^i(\theta)\hat{\zeta}^i\in \hat{{Q}}_{h,1},\\
\end{cases}
\end{equation}
where $\hat{\bm{\varphi}}^i$ and $\hat{\zeta^i}$ are the \emph{parameter--independent} basis functions of the FE spaces $\hat{V}_{h,2}^0$ and $\hat{Q}_{h,1}$ respectively.
We can now introduce the following forms:
\begin{eqnarray*}
\hat{a}_h(\hat{\bm{\varphi}}^i, \hat{\bm{\varphi}}^q; \theta)&:=&{\mu(\nabla\hat{\bm{\varphi}}^i, \nabla\hat{\bm{\varphi}}^q)_{\mathcal{R}} - \mu(\nabla\hat{\bm{\varphi}}^i\bm{n}_{\Gamma}, \hat{\bm{\varphi}}^q)_{\Gamma_{\theta}}}\nonumber\\
&&{-\mu(\hat{\bm{\varphi}}^i, \nabla\hat{\bm{\varphi}}^j\bm{n}_{\Gamma})_{\Gamma_{\theta}}
+ {\gamma\mu}/{h}(\hat{\bm{\varphi}}^i, \hat{\bm{\varphi}}^q)_{\Gamma_{\theta}}}
\\
&& + {\gamma\phi}/{h}( \hat{\bm{\varphi}}^i \cdot \bm{n}_\Gamma,\hat{\bm{\varphi}}^q\cdot \bm{n}_\Gamma)_{\Gamma_{\theta}}
\nonumber \\
\hat{c}_h(\hat{\bm{\varphi}}^k; \hat{\bm{\varphi}}^i, \hat{\bm{\varphi}}^q; \theta) & := & {((\hat{\bm{\varphi}}^k\cdot\nabla)\hat{\bm{\varphi}}^i, \hat{\bm{\varphi}}^q)_{\mathcal{R}}}\nonumber\\
\hat{b}_h(\hat{\zeta}^i, \hat{\bm{\varphi}}^q; \theta)&:=& -{(\hat{\zeta}^i, \nabla\cdot\hat{\bm{\varphi}}^q)_{\mathcal{R}} + (\hat{\zeta}^i\bm{n}_T,\hat{\bm{\varphi}}^q)_{\Gamma_{\theta}}}\nonumber\\
%
%
\end{eqnarray*}
We can therefore define the following matrices:
\begin{eqnarray*}
{\hat{\bm{A}}(\theta)_{iq}} &:=& {\hat{a}_h(\hat{\bm{\varphi}}^i, \hat{\bm{\varphi}}^q)+\hat{c}_h(\hat{\bm{\ell}}(\theta);  \hat{\bm{\varphi}}^i, \hat{\bm{\varphi}}^q)}
\nonumber\\
&&{+ \hat{c}_h(\hat{\bm{\varphi}}^i; \hat{\bm{\ell}}(\theta), \hat{\bm{\varphi}}^q) 
+ {{{g}^{GP}_{\mu}(\hat{\bm{\varphi}}^i, \hat{\bm{\varphi}}^q)}}}
\nonumber\\
{\hat{\bm{N}}(\underline{\textbf{u}}_h(\theta); \theta)_{iq}} &:=& {\sum_{k=1}^{N_u^h(\theta)}{u}_h^k(\theta)\hat{c}_h(\hat{\bm{\varphi}}^k; \hat{\bm{\varphi}}^i, \hat{\bm{\varphi}}^q)+}
\nonumber\\
&&+{{{{g}}_u^{GP}(\sum_{k=1}^{N_u^h}{u}_h^k(\theta)\hat{\bm{\varphi}}^k; \hat{\bm{\varphi}}^i, \hat{\bm{\varphi}}^q)}+ {{g}}_\beta^{GP}(\sum_{k=1}^{N_u^h(\theta)}{u}_h^k(\theta)\hat{\bm{\varphi}}^k; \hat{\bm{\varphi}}^i, \hat{\bm{\varphi}}^q)}
\nonumber\\
{\hat{\bm{B}}(\theta)_{iq}} &:=& {\hat{b}_h(\hat{\zeta}^i, \hat{\bm{\varphi}}^q),}
\nonumber\\
{\hat{\bm{C}}(\theta)_{iq}} &:=& {{g}^{GP}_p(\sum_{k=1}^{N_u^h(\theta)}{u}_h^k(\theta)\hat{\bm{\varphi}}^k; \hat{\zeta}^i, \hat{\zeta}^q)}
\end{eqnarray*}
%
%
%
%
%
Thanks to the matrices introduced, we can now state the algebraic formulation of the Navier--Stokes problem, after the snapshots extension and after the lifting of the inlet condition:
$$
R({U}_h(\theta), \theta):=
  \begin{bmatrix} {\hat{\bm{A}}}(\theta)
+  {{\hat{\bm{N}}}}({{\underline{\textbf{u}}}^0_h}(\theta); \theta) & {\hat{\bm{B^T}}}(\theta)
  \\ {\hat{\bm{B}}}(\theta) 
  & \hat{\bm{C}}(\theta)
  \end{bmatrix}
  \begin{bmatrix} {{\underline{\textbf{u}}}^0_h}(\theta) \\ {\underline{p}}_h(\theta) \end{bmatrix}
  -
  \begin{bmatrix} {{{\hat{\bm{F}}}}_1} (\theta) \\ {{{\hat{\bm{F}}}_2}}(\theta) \end{bmatrix}
  {{= \begin{bmatrix} 0 \\ 0 \end{bmatrix}}}
  ,
$$
where ${U}_h(\theta) = ({\underline{\textbf{u}}}^0_h(\theta), {\underline{p}}_h(\theta))$, ${(\hat{\bm{F}}_1(\theta))_i}:=\int_{\mathcal{R}}\hat{\bm{\varphi}}^i\cdot\bm{f}\,dx + \hat{a}_h(\hat{\bm{\varphi}}^i, \hat{\bm{\ell}}(\theta);\theta) + \hat{c}_h(\hat{\bm{\ell}}(\theta); \hat{\bm{\ell}}(\theta), \hat{\bm{\varphi}}^i)$ and ${\hat{\bm{F}}_2(\theta)}=\bm{0}$ since the embedded Dirichlet boundary data $g_D=0$.
}
\subsubsection{Reduced Basis generation}\label{reduced basis generation}
Let us now denote by  $\theta^{(j)}$ each parameter
in a finite dimensional training set $\mathcal{P}_{train} = \{ \theta^1, \dots, \theta^M\}$ for a large number $M$.
The snapshots matrices $\bm{\mathcal{S}_u}$ and $\bm{\mathcal{S}_p}$ for the fluid velocity and the fluid pressure are defined as follows:
\begin{gather}
 \bm{\mathcal{S}_u}= [\hat{\bm{u}}^0_h(\theta^1),\dots,\hat{\bm{u}}^0_h(\theta^M)] \in \mathbb{R}^{N_u^h\times M},\quad\bm{\mathcal{S}_p} = [\hat{{p}}_h(\theta^1),\dots,\hat{{p}}_h(\theta^M)] \in \mathbb{R}^{N_p^h\times M}.
\end{gather}

In order to make the pressure approximation stable at the reduced order level we also introduce a velocity supremizer variable $\bm{s_h}$: see \cite{ballarin2015supremizer,RoHuMa13,RoVe07} for a more detailed introduction to the supremizer enrichment for Navier--Stokes equation. 
{
We start with a Laplace problem for the supremizer $\bm{s}(\theta)$, $\forall \theta\in\mathcal{P}$:
\begin{equation}\label{Poisson supremizer}
\begin{cases}
-\Delta\bm{s}(\theta)=-\nabla p_h(\theta)\text{ in }\Omega(\theta),\\
\bm{s}(\theta)=0\text{ in }\Gamma_{in\cup D\cup out},\\
\bm{s}(\theta)=0\text{ in }\Gamma_{\theta}
\end{cases}
\end{equation}
Again, we impose boundary conditions \eqref{Poisson supremizer}$_2$ strongly, and \eqref{Poisson supremizer}$_3$ weakly, thanks to the Nitsche method. We have therefore the following discretized problem: find $\bm{s}_h(\theta)\in V_{h,k}^0(\theta)$ such that $\forall \bm{v}_h(\theta)\in V_{h, k}^0(\theta)$:
\begin{eqnarray*}
(\nabla\bm{s_h}(\theta), \nabla\bm{v_h}(\theta))_{\Omega(\theta)} - (\nabla\bm{s_h}(\theta)\bm{n}_{\Gamma}, \bm{v_h}(\theta))_{\Gamma_{\theta}} - (\nabla\bm{v_h}(\theta)\bm{n}_{\Gamma}, \bm{s_h}(\theta))_{\Gamma_{\theta}} + \frac{\lambda_s}{h}(s_h(\theta), v_h(\theta))_{\Gamma_{\theta}} \\
+ \mathit{g}^{GP}
 (\bm{s_h}(\theta), \bm{v_h}(\theta))
= -(\nabla p_h(\theta), \bm{v}_h(\theta))_{\Gamma_{\theta}},
\end{eqnarray*}
where the Ghost Penalty term is given by:
\begin{equation*}
\mathit{g}^{GP}(\bm{s_h}(\theta), \bm{v_h}(\theta)) = \sum_{1\leq j\leq 2}\sum_{F\in\mathcal{F}_{\Gamma_{\theta}}}\gamma_s^jh^{2j +1}(\ldb\partial_n^j \bm{s_h}(\theta)\rdb, \ldb\partial_n^j\bm{v_h}(\theta)\rdb)_F,
\end{equation*} }
We employ the same natural smooth extension (and the same extended FE space used for velocity) also for the supremizer, thus obtaining the extendend snapshots $\hat{\bm{s}}_h$.
These snapshots are then collected in the snapshot matrix
\begin{gather*}
{\bm{\mathcal{S}}_s} = [\hat{\bm{s}}_h(\theta^1),\dots,\hat{{\bm{s}}}_h(\theta^M)] \in \mathbb{R}^{N_u^h\times M},
\end{gather*}

We then carry out a compression by POD on the snapshots matrices, namely $\bm{\mathcal{S}}_u$, $\bm{\mathcal{S}}_s$ and $\bm{\mathcal{S}}_p$, following e.g. \cite{Kunisch2002492}.
This derives an eigenvalue problem, that for the velocity for example reads:
\begin{gather}
{\bm {\mathcal{C}}^u}\bm{Q}^u = \bm{Q}^u\bm{\Lambda}^u ,\quad \mbox{\hspace{0.5cm} for }{\bm{\mathcal{C}}}^u_{ij} = ({\hat{\bm{u}}_h(\theta^i),\hat{\bm{u}}_h(\theta^j)}) _{L^2{(\hat{\mathcal{I}}_h )}} \mbox{,\, } i,j = 1,\dots,M ,\nonumber
\end{gather}
where $\bm{\mathcal{C}^u}$ is the correlation matrix derived from the $\theta$-independent snapshots, $\bm{Q}^u$ is an eigenvectors square matrix and $\bm{\Lambda}^u$ is a diagonal matrix of eigenvalues. Similar eigenvalue problems can be derived for the supremizer and for the pressure. {We refer the reader interested in the details about the POD and its implementation to \cite{HesthavenRozzaStamm, ballarin2015supremizer}. The $i$--th reduced basis function $\bm{\Phi}_i^u$ for the fluid velocity, for example, is then obtained (possibly after $L^2$ normalization) by applying the snapshots matrix $\bm{\mathcal{S}}_u$ to the $i$--th column of the matrix $\bm{Q}^u$.}
\subsection{Online algebraic system}\label{sec:online}
Thanks to the POD on the velocity snapshots and thanks to the enrichment with supremizer snapshots we obtain a set $\{\bm{\Phi}_1^u, \dots, \bm{\Phi}_N^u,\bm{\Phi}_1^s, \dots, \bm{\Phi}_N^s\}$ of $2N$ basis functions for the reduced order approximation of the velocity, and a set $\{\Phi_1^p, \dots, \Phi_N^p\}$ of $N$ basis functions for the reduced order approximation of the pressure.
We define: $\hat{V}_N$, the enriched reduced basis space for the velocity, and $\hat{Q}_N$, the reduced basis space for the pressure:
\begin{align*}
\hat{V}_N = \text{span}\{\bm{\Phi}_1^{u, s}, \hdots, \bm{\Phi}_{2N}^{u,s}\}, \qquad
\hat{Q}_N = \text{span}\{\Phi_1^p, \hdots, \Phi_N^p\}
\end{align*}
where $N< M$ is chosen according to the eigenvalue decay of $\bm{\Lambda}_{ii}^u$ and $\bm{\Lambda}_{ii}^p$, see for instance \cite{Rozza2008229,BeOhPaRoUr17}. In the definition of the reduced space for the fluid velocity, we used a unified notation: $\bm{\Phi}_i^{u,s}=\bm{\Phi}_i^u$ for $i=1,\dots, N$ and  $\bm{\Phi}_i^{u,s}=\bm{\Phi}_i^s$ for $i=N+1,\dots, 2N$. We also remark the fact that the finite dimensional reduced spaces $\hat{V}_N$ and $\hat{Q}_N $ are \emph{parameter--independent}, thanks to the natural smooth extension performed on the snapshots.
We can now introduce the online velocity $\bm{u}_N(\theta)$ and the online pressure $p_N(\theta)$:
\begin{gather}
\bm{u}_N(\theta) := \sum_{i=1}^{2N}\underline{\textbf{u}}_N^i(\theta)\bm{\Phi}_i^{u, s} = \bm{L}_{u,s}\underline{\textbf{u}}_N(\theta), \label{reduced velocity}\\
p_N(\theta) := \sum_{i=1}^N\underline{p}^i_N(\theta)\Phi_i^p = \bm{L}_p\underline{p}_N(\theta), \label{reduced pressure}
\end{gather}
where $\bm{L}_{u,s}\in\mathbb{R}^{N_u^h\times 2N}$ and $\bm{L}_p\in\mathbb{R}^{N_p^h\times N}$ are rectangular matrices containing the FE degrees of freedom of the basis of $\hat{V}_N$ and $\hat{Q}_N$.
The parameter {\textit{dependent}} solution vector $\underline{\textbf{u}}_N(\theta)\in\mathbb{R}^{2N}$ and $\underline{p}_N(\theta) \in \mathbb{R}^{N}$ and the parameter {\textit{independent}} reduced basis functions $\bm{\Phi}_i^{u, s}$, $\Phi_i^p$ are the key ingredients necessary to
perform a Galerkin projection of the full system onto the aforementioned reduced basis space.
By introducing the vector $U_N(\theta) = (\underline{\textbf{u}}_N(\theta), \underline{p}_N(\theta))$,  the algebraic formulation, at the reduced order level, of the steady Navier--Stokes problem, reads as follows:
\begin{eqnarray}\label{online algebraic system}
\hat{R}(U_N(\theta),\theta):=
 \begin{bmatrix} {\bm{L}_{u, s}^T(\hat{\bm{A}}}(\theta) + {\hat{\bm{N}}(\underline{\textbf{u}}_N(\theta); \theta))\bm{L}_{u, s}} 
 & {\bm{L}_{u, s}^T\hat{\bm{B}}}^T(\theta)\bm{L}_p
  \\ {\bm{L}_p^T\hat{\bm{B}}}(\theta)\bm{L}_{u, s}
  & \bm{L}_p^T\hat{\bm{C}}(\theta)\bm{L}_p
  \end{bmatrix}
  U_N(\theta) -
  \begin{bmatrix}
   \bm{L}_{u, s}^T \hat{\bm{F}}_1(\theta)\\
  \bm{L}_p^T\hat{\bm{F}}_2(\theta)
  \end{bmatrix}=0.
\end{eqnarray}

{We point out that the matrices appearing in \eqref{online algebraic system} have huge kernels, indeed, as we can see, we are considering the FE basis functions of the FE spaces $\hat{V}_{h, 2}$ and $\hat{Q}_{h,1}$ that are parameter--independent, and defined on the whole background mesh.
The aforementioned formulation \eqref{online algebraic system} however is required only by the ROM procedure: in particular, during the solution of the online reduced system \eqref{online algebraic system}, we are going to discard all the  entries in the matrices that are associated with DOFs that are situated outside of the computational domain $\Omega_h^*(\theta)$. This means that the value of the reduced solution $\bm{u}_N(\theta)$ outside $\Omega_h^*(\theta)$ is not interesting and can be discarded during the analysis of the numerical results.}

In the above POD-ROM solution, we clarify that we have to assemble the matrices of the high fidelity system. For a ``cheaper'' in time execution and less computation resources costs, one could achieve further improvement employing  hyper reduction techniques as in \cite{Xiao20141,BARRAULT2004667,Carlberg2013623,stabile_geo_}.
\subsection{Numerical results}\label{numerical results steady}
In this paragraph we present the results obtained by applying the aforementioned reduction techniques to our model problem. For our simulation, the fluid viscosity is $\mu_f = 0.05$ $cm^2/s$ and the fluid density is $
1$ $g/cm^3$. The prescribed inlet velocity is given by $\bm{u}_{in}=(1, 0)$ $m/s$. The reduced basis have been obtained with a Proper Orthogonal Decomposition on the set of snapshots: this reduction technique, although costly in computational terms, is very useful as it gives an insight on the rate of decay of the eigenvalues related to each component of the solution.
We take $N_{\text{train}} = 150$, and we generate randomly $N_{\text{train}}$ uniformly distributed values for the parameter $\theta$. We then run a POD on the collected set of snapshots and we obtain our basis functions, with which we are going to compute the reduced solutions $(\bm{u}_N(\theta^{(i)}), p_N(\theta^{(i)}))$, where $i=1, \dots, N_{\text{test}}$, and $N$ is the number of basis functions that we use.
\begin{figure}
\begin{center}
\begin{tabular}{cc}
\subfigure[First six modes for the velocity.]{
\includegraphics[scale=0.2]{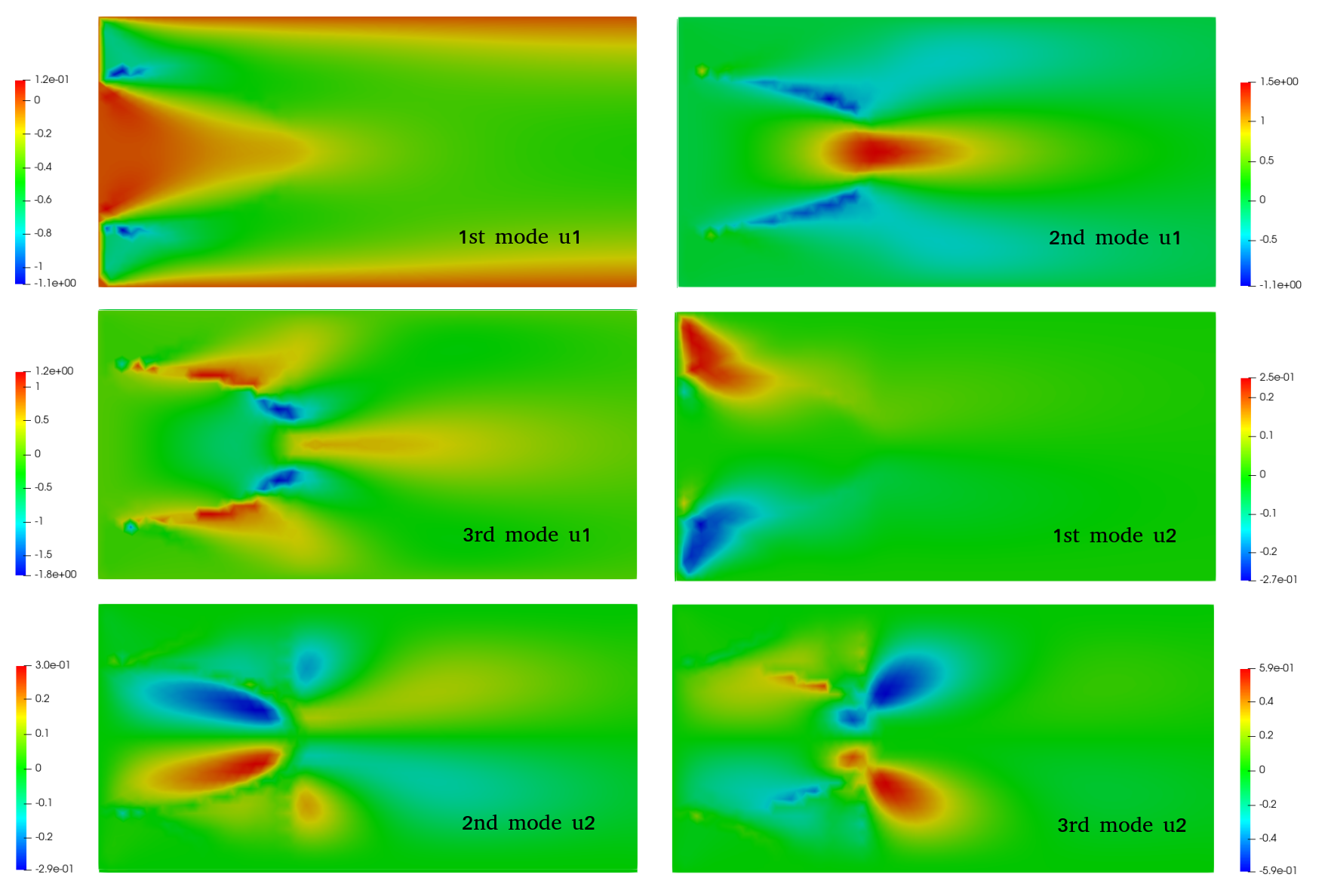}}
\label{velocity modes}
\\
\subfigure[First four modes for the pressure.]{
\includegraphics[scale=0.2]{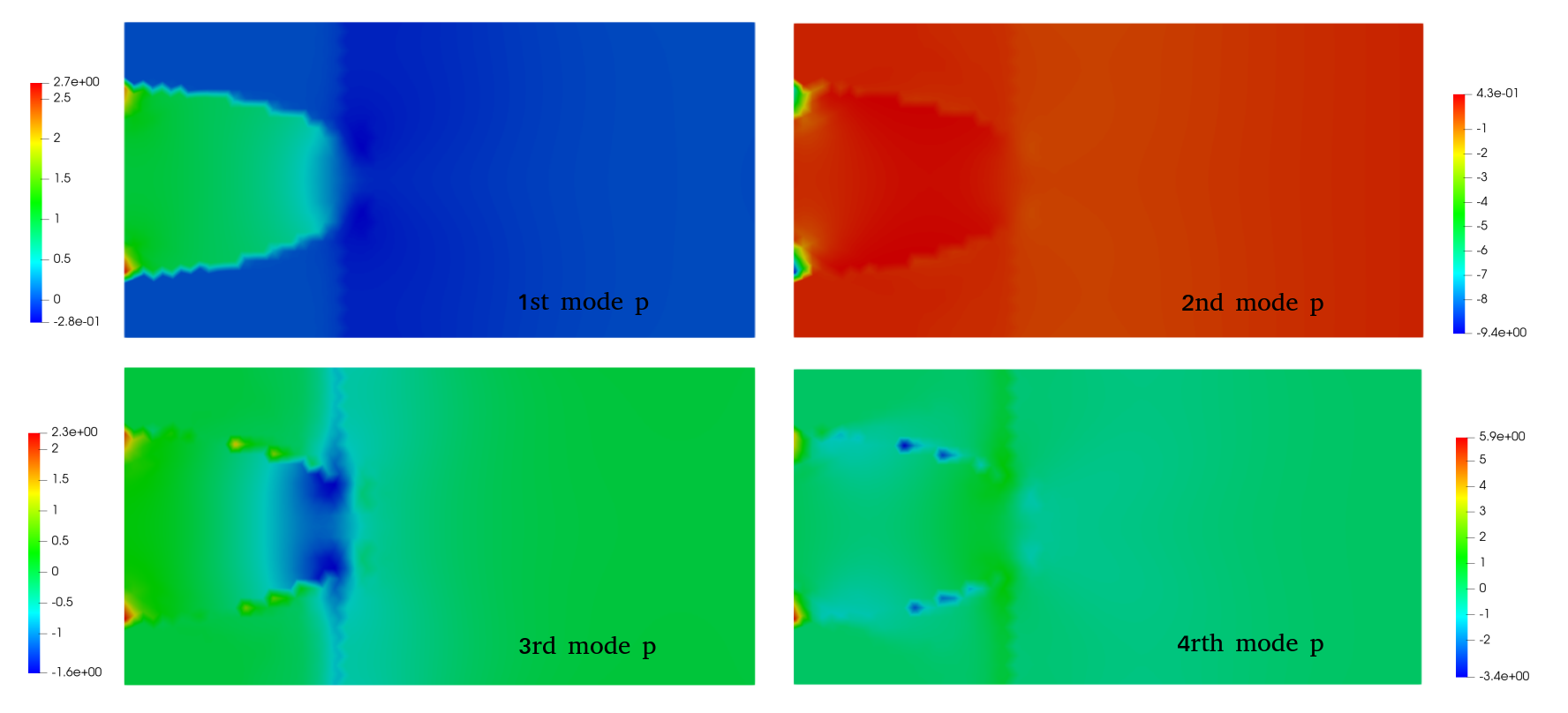}}
\label{pressure modes}
\end{tabular}
\caption{Steady system: Some reduced basis modes for velocity and pressure for a geometrically patrametrized  Navier--Stokes system.}\label{complex_fig}
\end{center}
\end{figure}
\begin{figure}
\centering
\includegraphics[scale=0.5]{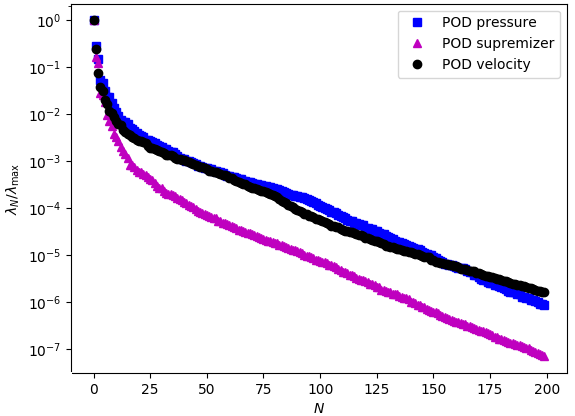}
\caption{Steady case: POD eigenvalues decay for the fluid velocity $\textbf{u}$ (black), the fluid pressure $p$ (blue), and the fluid supremizer $\textbf{s}$ (magenta), for a set of $N_{train}=200$ snapshots.}
\label{POD eigenvalues}
\end{figure}
\begin{figure}
\centering
\includegraphics[scale=0.55]{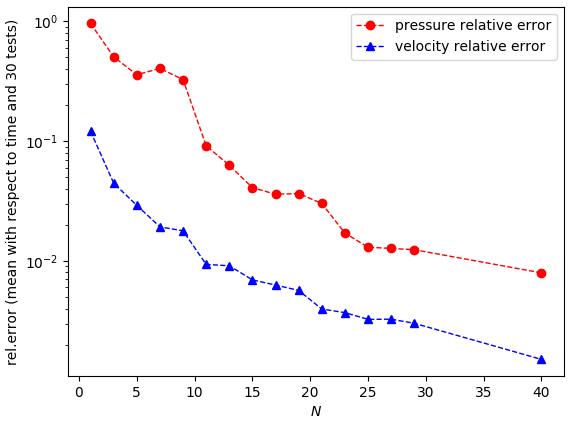}
\caption{Steady case: Mode dependent errors between high fidelity and reduced order approximation, with the supremizer enrichment.}
\label{r_err_no_supremizer}
\end{figure}
Figures \ref{complex_fig}(a) and \ref{complex_fig}(b) give an example of the first modes that we obtain with this procedure, whereas in Figure \ref{POD eigenvalues} we report the decay of the eigenvalues for all the components of the solution and for the supremizer.
To test the reduced order model we generate randomly $N_{\text{test}}=30$ uniformly distributed values values for $\theta\in\mathcal{P}_{\text{test}}$. We are interested in the behavior of the relative approximation error that we obtain by changing the number of basis functions $N$ used to build the reduced solution.
 In order to do this we let $N$ vary in a discrete set $\mathcal{N}$: for a fixed value of $N\in\mathcal{N}$, and for each $\theta^i$, $i=1, \dots, N_{\text{test}}$, 
we compute both the reduced solution $(\bm{u}_N(\theta^i), p_N(\theta^i))$ and the corresponding full order solution $(\bm{u}_h(\theta^i), p_h(\theta^i))$. We compute the $L^2$ relative error $\epsilon_u^{N, i}$ for the velocity and the relative error $\epsilon_p^{N, i}$ for the pressure; then we compute the average approximation errors $\overline{\epsilon}_u^N$ and $\overline{\epsilon}_p^N$ for every $N\in\mathcal{N}$, defined as:
\begin{equation*}
\overline{\epsilon}_u^N = \frac{1}{N_{\text{test}}}\sum_{i=1}^{N_{\text{test}}}\epsilon_u^{N, i}.
\end{equation*}
Figure \ref{r_err_no_supremizer} shows the relative approximation errors plotted against the number $N$ of basis functions used, with the use of the supremizer enrichment at the reduced order level.
\begin{figure} 
\begin{center}
\includegraphics[scale=0.25]{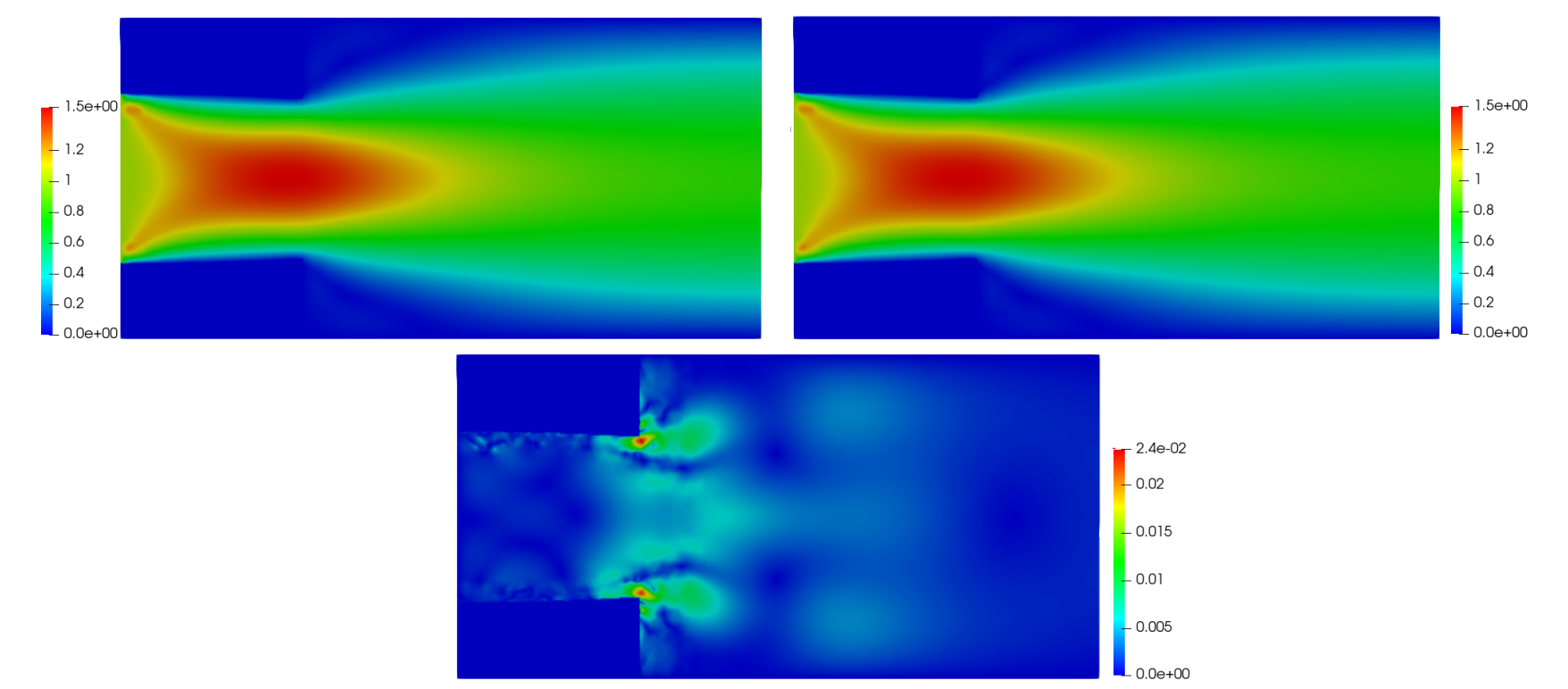}
\caption{Steady case: Uncut geometry and the high fidelity  velocity solution for parameter $\theta=-0.015854$ (left), reduced order solution for the same $\theta$ (right) and approximation error (middle)}\label{velocity_approx_error_with_s_}
\end{center}
\end{figure}
\begin{figure}
\begin{center}
\includegraphics[scale=0.25]{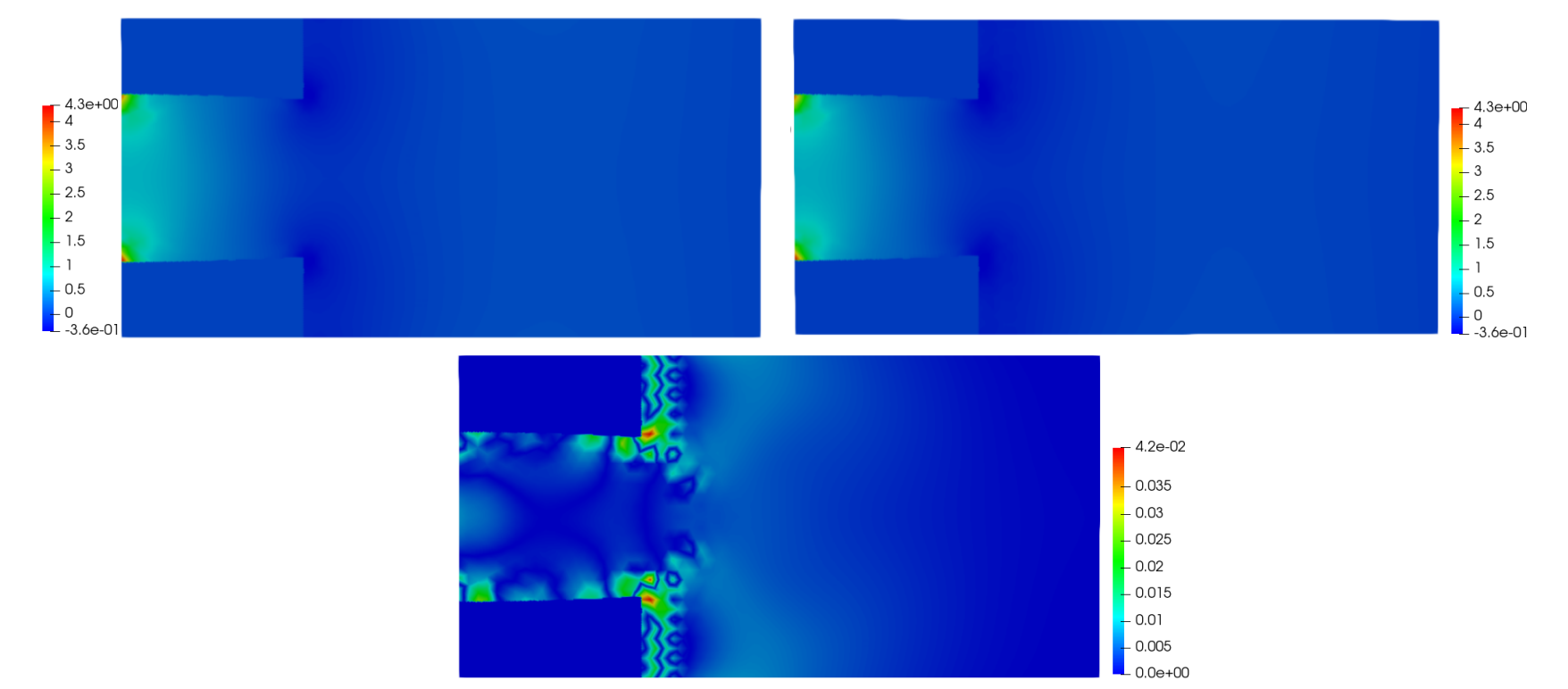}
\caption{Steady case: Uncut geometry and the high fidelity pressure solution for parameter $\theta=-0.015854$ (left), reduced order solution for the same $\theta$ (right) and approximation error (middle)}\label{pressure_approx_error_with_s}
\end{center}
\end{figure}

Figures \ref{velocity_approx_error_with_s_} and 
 \ref{pressure_approx_error_with_s} show the approximation error for the velocity and pressure, for a given test value of the parameter, with the supremizer enrichment. 
  It is worth to mention that the approximation error tends to concentrate near the cut between the physical domain and the background mesh, similar to to experiments in the works of
\cite{KaBaRO18,KaratzasStabileNouveauScovazziRozza2018,KaratzasStabileNouveauScovazziRozzaNS2019,KaratzasRozzaCH2019}, phenomenon which will be studied in a future work.
{
\subsubsection{Integration over boundary elements}
The simulations presented in this manuscript have been implemented with \emph{ngsxfem}, which is an add-on library to the finite element package Netgen/NGSolve, which enables the use of unfitted finite element technologies, such as CutFEM indeed. One of the main tools of ngsxfem is the availability of a routine that "studies" the topology of the domain, identifying elements of the background mesh that are situated inside or outside the levelset function, or are cut-elements. Thanks to this routine, ngsxfem then takes automatically care of the integration of bilinear forms over cut-elements.
}

\section{Unsteady Navier--Stokes}\label{Unsteady NS}
We now extend the previous treatment to the unsteady incompressible Navier--Stokes problem, by introducing the time evolution term $\partial_t{\bm{u}}(\theta)$ in system (\ref{eq:NS_strong}). The physical domain over which the problem is formulated is always the one in Figure \ref{mesh}.
\subsection{Strong formulation}\label{strong formulation unsteady}
Given a time interval of interest $[0, T]$, the strong formulation of the problem reads as follows: for every $\theta\in\mathcal{P}$, and for every $t\in[0, T]$, find $\bm{u}(t;\theta)\colon\Omega(\theta)\mapsto\mathbb{R}^2$ and $p(t;\theta)\colon\Omega(\theta)\mapsto\mathbb{R}$ such that:
\begin{equation}\label{eq:NS_time_strong}
\begin{cases}
\partial_t{\bm{u}}(\theta)-\mu\Delta \bm{u}(\theta) + \nabla p(\theta) + (\bm{u}(\theta)\cdot \nabla)\bm{u}(\theta) =  \bm{f}(\theta) \quad \text{in $\Omega(\theta)\times[0,T]$}, \\
\text{div}\bm{u}(\theta)=0 \quad \text{in $\Omega(\theta)\times[0,T]$},\\
\bm{u}(\bm{x},0;\theta)={\bm{u}}^0(\bm{x}, \theta) \quad \text{in $\Omega(\theta)$},
\end{cases}
\end{equation}
with geometrical parameterization identical to that in the previous subsection. 
{The problem is completed by the following boundary conditions: a prescribed inlet velocity $\bm{u}_{in}$ at the inlet boundary $\Gamma_{in}(\theta)$, a zero outflow condition at the outflow boundary, a boundary condition $\bm{u}(\theta)\cdot\bm{n}_D=0$ on $\Gamma_D$ and a no--slip boundary condition onf $\Gamma_{\theta}$.}

 \subsection{Space discretization and time--stepping scheme}\label{Time discretization}
{We discretize in time by a backwards Euler approach: we discretize the time interval $[0, T]$ with the following partition:}
\begin{equation*}
{0=t^0 <  \ldots < t^{N_t}=T},
\end{equation*}
{where every interval $(t^{n},t^{n+1}]$ has measure $\tau_{n+1}=t^{n+1}-t^{n}$, $n=0,...,N_t-1$.
The discrete version of the initial condition $\bm{u}^0(x; \theta)$ is denoted by $\bm{u}^0_h(x;\theta)$; we denote $\bm{u}_h^n$ the discrete fluid velocity at time step $t^n$, and similar notation is used for the pressure.}

After having applied a time stepping scheme, the space discretized weak formulation of the problem reads as follows: for every $n=0, \hdots, N_t-1$, we seek a discrete velocity ${\bm{u}}_h^{n+1}(\theta)\in V_{h,2}^D(\theta)$ and discrete pressure $p_h^{n+1}(\theta) \in Q_{h,1}(\theta)$, such that for every $(\bm{v}_h(\theta), q_h(\theta))\in V_{h,2}^0(\theta)\times Q_{h,1}(\theta)$, it holds:
\begin{equation*}
m(\bm{u}_h^{n+1}(\theta) - \bm{u}_h^n(\theta), v_h) +
\tau_{n+1} \mathcal{A}(\bm{u}_h^{n+1}(\theta),p_h^{n+1}({\theta}), \bm{v}_h(\theta),q_h(\theta)) =
\tau_{n+1} \mathcal{L}(\bm{v}_h(\theta)),
\end{equation*}
where
\begin{equation*}
m(\bm{w}_h, \bm{v}_h):=\int_{\Omega(\theta)}\bm{w}_h\cdot\bm{v}_h\,dx,
\end{equation*}
 where $\mathcal{L}(\bm{v}_h(\theta))$ is defined as in Section \ref{discrete steady Stokes weak formulation}, and $\mathcal{A}({\bm{u_h}}^{n+1}(\theta),p_h^{n+1}(\theta), \bm{v_h}(\theta), q_h(\theta))$  is defined as follows:
{
\begin{eqnarray*}
\mathcal{A}(\bm{u}_h(\theta),p_h(\theta); \bm{v}_h(\theta), q_h(\theta)) &=& {a_h(\bm{u}_h(\theta), \bm{v}_h(\theta)) + b_h(p_h(\theta), \bm{v}_h(\theta)) - b_h(q_h(\theta), \bm{u}_h(\theta)) +}\nonumber\\
&&{+ c_h(\bm{u}_h(\theta); \bm{u}_h(\theta), \bm{v}_h(\theta) + {g}_u^{GP}(\bm{u}_h(\theta); \bm{u}_h(\theta), \bm{v}_h(\theta)) + g_{\mu}^{GP}(\bm{u}_h(\theta), \bm{v}_h(\theta)) + }\nonumber\\
&&{g_{\sigma}^{GP}(\bm{u}_h(\theta), \bm{v}_h(\theta)) + g_p^{GP}(\bm{u}_h(\theta); p_h(\theta), q_h(\theta))}
\end{eqnarray*}
The forms $a_h$, $b_h$, $c_h$ as well as the stabilization terms ${g}_u^{GP}$, $g_{\mu}^{GP}$, $g_{\beta}^{GP}$ and $g_p^{GP}$ have been introduced in Section \ref{discrete steady Stokes weak formulation}. 
}

\subsection{POD and reduced basis generation}\label{POD unsteady}
Similarly to what has been done in the previous Section \ref{sec:ROM}, in the time dependent case an \textit{offline/online} procedure will be employed, that will lead to the generation of a proper reduced basis set.
Since the system is both  (geometrical) parameter and time-dependent, we sample not only the geometrical parameter $\theta$, but also the time $t$, with the sample points $t^k \in \{t^0,\dots,t^{N_t}\} \subset [0,T]$.
This procedure is computationally more expensive and results in a much larger total number of snapshots to be collected with respect to the static system: the total number of snapshots that we collect is now equal to $\widehat{M} = M\cdot N_t$.
The snapshots matrices $\bm{\mathcal{S}_u}$, $\bm{\mathcal{S}_s}$ and $\bm{\mathcal{S}}_p$ are then given by:
\begin{gather}
\bm{\mathcal{S}}_u = [\hat{{\bm{u}}}^0_h(\theta^1,t^0),\dots,\hat{{\bm{u}}}^0_h(\theta^1,t^{N_t}),\dots,\hat{{\bm{u}}}^0_h(\theta^M, t^0), \dots, \hat{{\bm{u}}}^0_h(\theta^M,t^{N_t})] \in \mathbb{R}^{N_u^h\times \widehat{M}},\\
\bm{\mathcal{S}}_s = [\hat{{\bm{s}}}_h(\theta^1,t^0), \dots, \hat{{\bm{s}}}_h(\theta^1,t^{N_t}),\dots, \hat{{\bm{s}}}_h(\theta^M, t^0), \dots, \hat{{\bm{s}}}_h(\theta^M,t^{N_t})] \in \mathbb{R}^{N_u^h\times \widehat{M}},\\
\bm{\mathcal{S}}_p = [\hat{{p}}_h(\theta^{(1)},t^0),\dots, \hat{{p}}_h(\theta^1,t^{N_t}), \dots, \hat{{p}}_h(\theta^M, t^0), \dots, \hat{{p}}_h(\theta^M,t^{N_t})] \in \mathbb{R}^{N_p^h\times \widehat{M}},
\end{gather}
{where we used the $\hat{}$ notation to indicate that we implemented, also in this case, a natural smooth extension of the FE solutions obtained in the offline phase of the method, in order to work with parameter--independent reduced basis functions and parameter--independent reduced basis spaces, as well as a lifting function for the inlet velocity, see Section \ref{lifting function}.}
We solve an eigenvalue problem like the one introduced in Section \ref{sec:ROM}, and finally, adopting
the notation of Section \ref{sec:online} we end up with the parameter--independent reduced basis spaces
\begin{align*}
\hat{V}_N = \text{span}\{\bm{\Phi}_1^{u,s}, \hdots, \bm{\Phi}_{2N}^{u,s}\}, \qquad
\hat{Q}_N = \text{span}\{\Phi_1^p, \hdots, \Phi_N^p\},
\end{align*}
Let us now denote by $(\bm{u}_N^n(\theta), p_N^n(\theta))$ the reduced solution at time-step $t^n$, for $n=0, \hdots, N_t$, where $\bm{u}_N^n(\theta):=\bm{u}_N^{n, 0}(\theta)+\ell(\theta)$, and $\bm{u}_N^{n, 0}(\theta)$ and $p_N^n(\theta)$ are defined as in \eqref{reduced velocity} and in \eqref{reduced pressure}, respectively. We can derive the subsequent reduced algebraic system for the unknown $U_h^{n+1}(\theta) = (\bm{u}_N^{n+1,0}(\theta), p_N^{n+1}(\theta))$:
{
\begin{equation}
\label{eq:system_linear_reduced}
\begin{bmatrix}
{\bm{L}_{u,s}^T\hat{\bm{M}}\bm{L}_{u, s}} & \bm{0}\\
\bm{0} & \bm{0}
\end{bmatrix}
U_N^{n+1}(\theta)
+
\tau_{n+1}
\hat{R}_N(U_N^{n+1}(\theta); \theta)
=
\begin{bmatrix}
{{\bm{L}_{u,s}^T\hat{\bm{M}}\bm{L}_{u, s}}} & \bm{0}\\
\bm{0} & \bm{0}
\end{bmatrix}
U_N^n(\theta),
\end{equation}
where $\hat{M}_{ij}:=(\hat{\bm{\varphi}}^i,\hat{\bm{\varphi}}^j)_{\mathcal{R}}$. Here $\{\hat{\bm{\varphi}}^i\}_{i=1}^{N_u^h}$ are the FE basis functions of the FE parameter--independent spaces $\hat{V}_{h, 2}^D$ and $\hat{Q}_{h,1}$, as defined in Section \ref{natural smooth extension}. Here, $\hat{R}(U_N(\theta);\theta)$ is defined as in Section \ref{sec:online}.
}

\subsection{Numerical results}\label{numer_examples_unsteady}
\label{sec:Stationary_in_time_geometry}
In this paragraph we present the results obtained by applying the proposed reduction technique to a time dependent case. The time-step used in our simulation is $\tau = 0.011s$, and the final time is $T=0.7s$. The fluid viscosity is $\mu_f=0.05$ $cm^2/s$ and the fluid density is $
1$ $g/cm^3$. We impose a constant inlet velocity $\bm{u}_{\text{in}} = (1, 0)$ $m/s$. We now take $N_{\text{train}} = 200$, and we generate randomly $N_{\text{train}}$ uniformly distributed values for the parameter $\theta$. We also remind that we sample the time interval $[0, T]$ with an equispaced sampling $\{t_0, \hdots, t_{N_t}\}$.
We then run a POD on the set of snapshots collected, and we obtain our basis functions with which we are going to compute the reduced solutions $(\bm{u}_N(t, \theta^i), p_N(t, \theta^i)$, where $i=1, \dots, N_{\text{test}}$, and $N$ is the number of basis functions that we use.
\begin{figure}
\begin{center}
\begin{tabular}{cc}
\subfigure[First six modes for the velocity.]{
\includegraphics[scale=0.2]{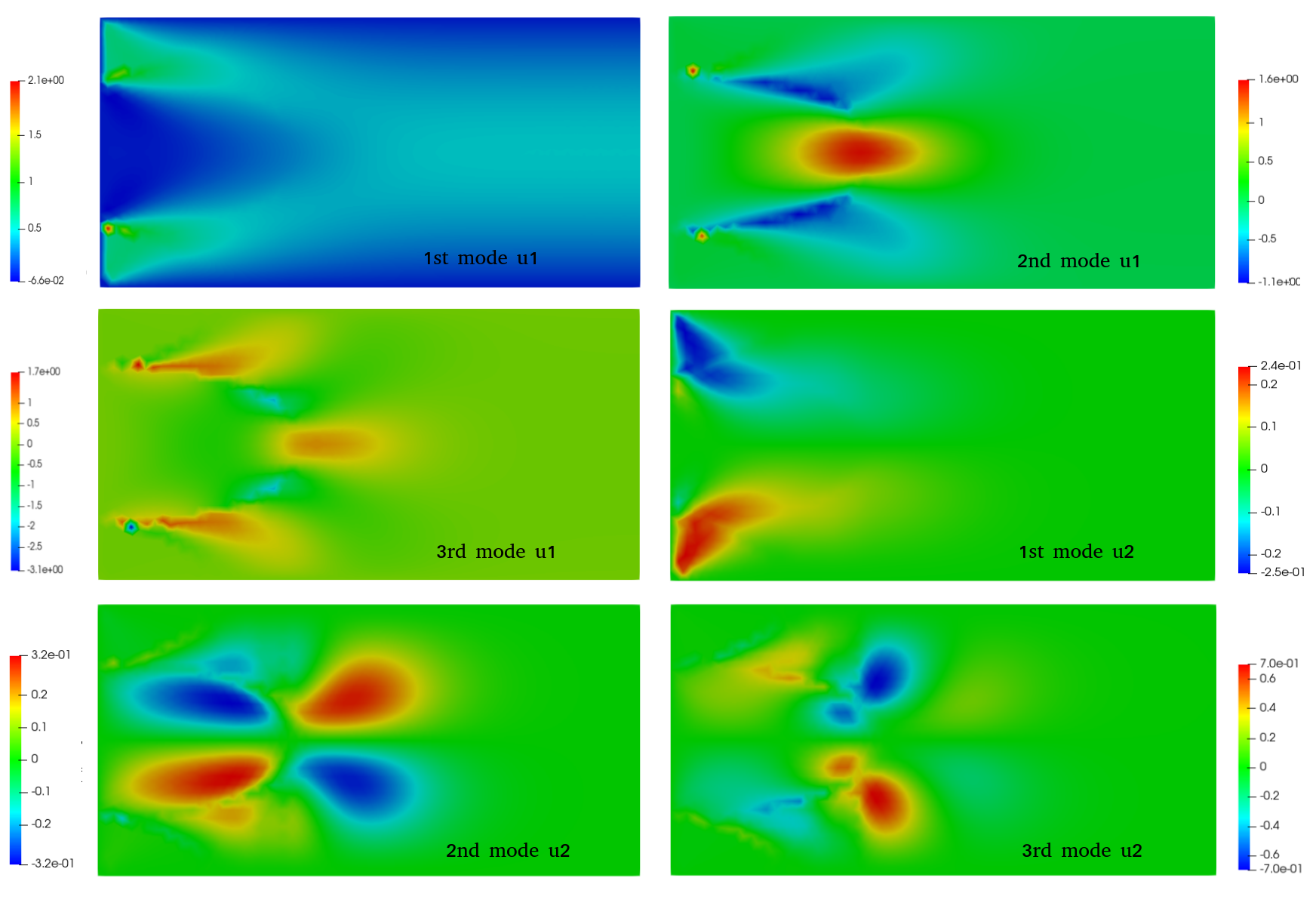}}
\\
\subfigure[First four modes for the pressure.]{
\includegraphics[scale=0.2]{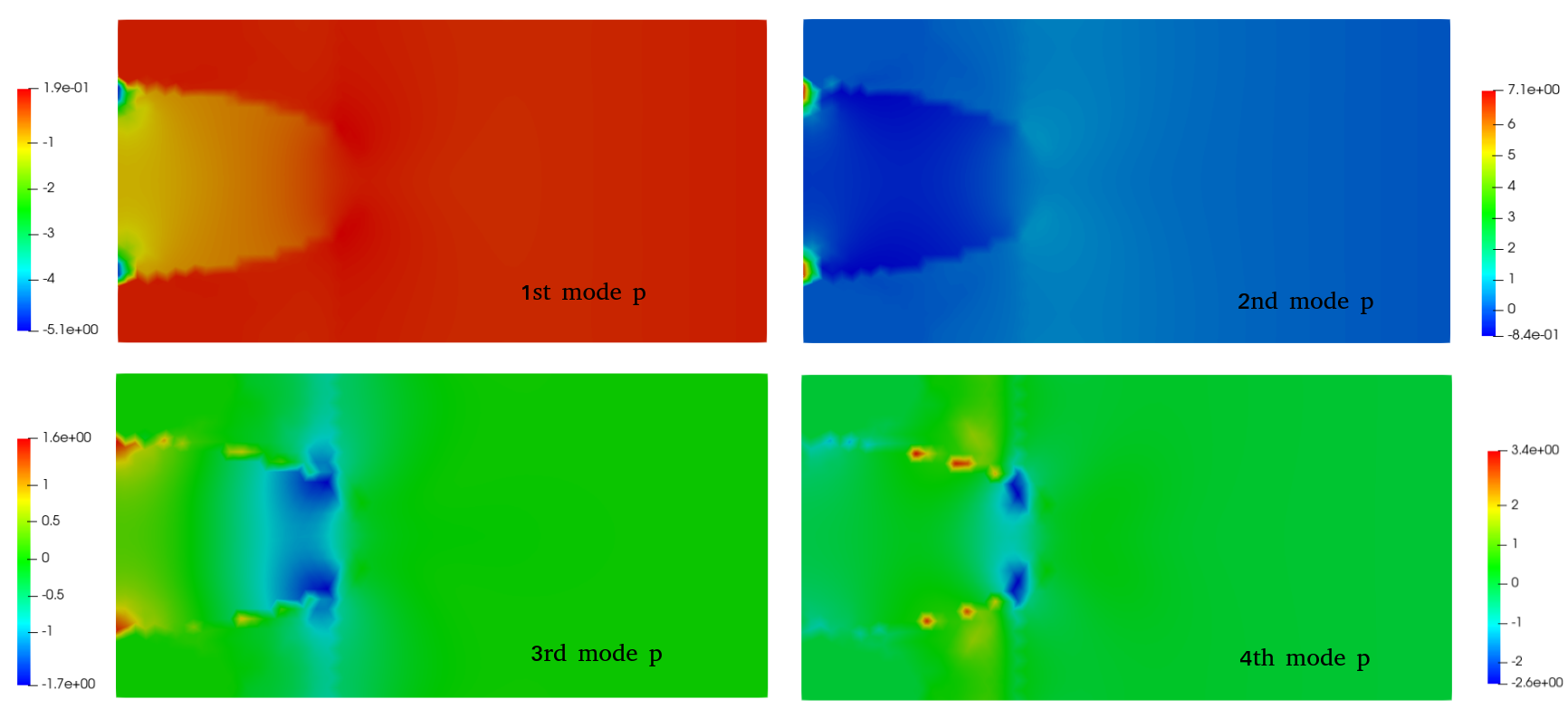}}
\end{tabular}
\caption{Unsteady system with static in time geometry: Some reduced basis modes for velocity and pressure for the evolutionary in time, geometrically patrametrized Navier--Stokes system.}
\label{t_velocity modes}
\label{t_pressure modes}
\end{center}
\end{figure}

\begin{figure}
\centering
\includegraphics[scale=0.5]{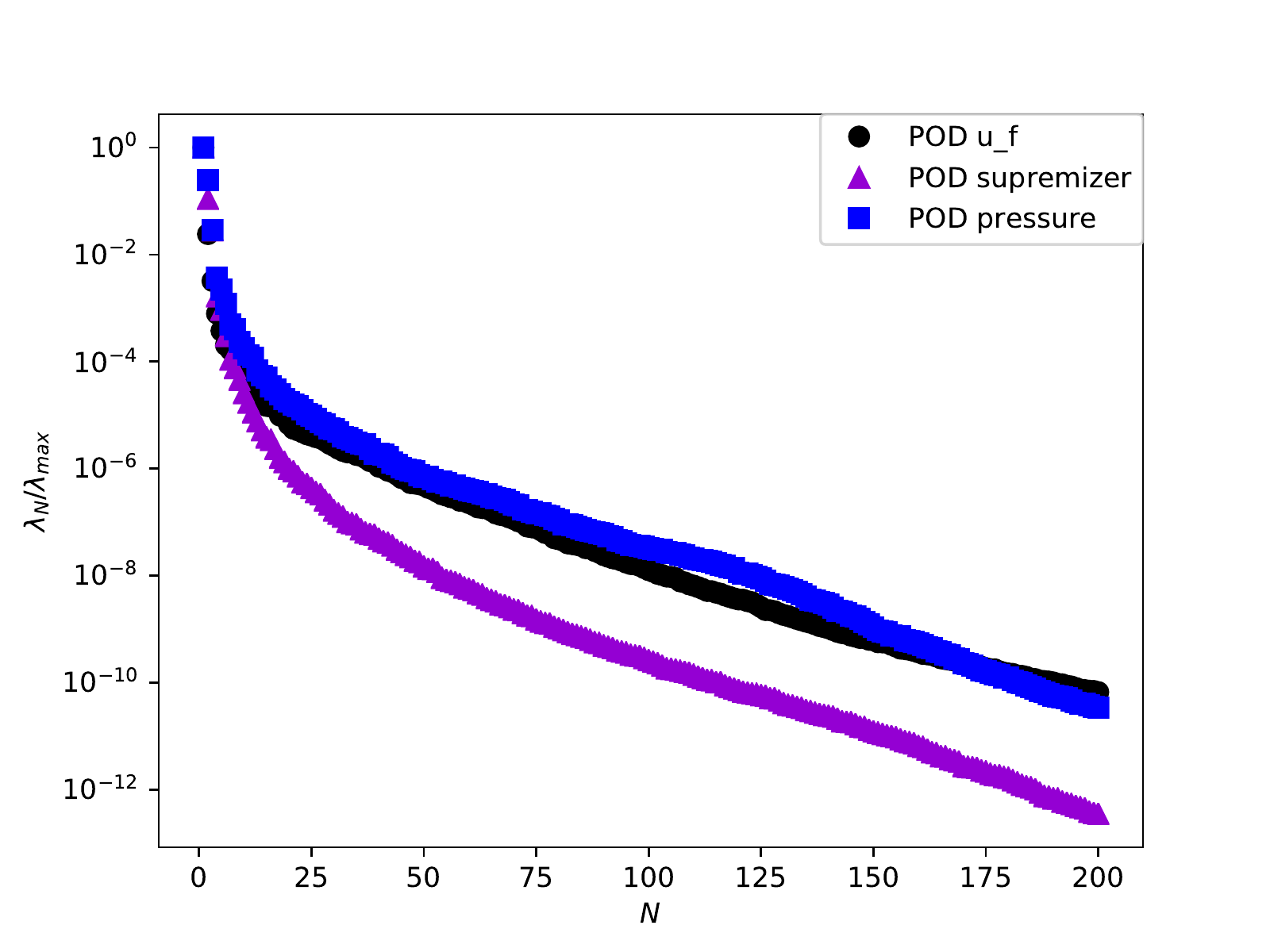}
\caption{Time dependent case and the POD eigenvalues decay for the fluid velocity $\textbf{u}$ (black), the fluid pressure $p$ (blue), and the fluid supremizer $\textbf{s}$ (magenta), for a set of $N_{train}=200$ snapshots.}
\label{t_POD eigenvalues}
\end{figure}

Figure \ref{t_velocity modes} gives an example of the first modes that we obtain with this procedure, whereas in Figure \ref{t_POD eigenvalues} we show the rate of decay of the eigenvalues for all the components of the solution and for the supremizer.
To test the reduced order model we generate randomly $N_{\text{test}}=30$ uniformly distributed values for $\theta\in\mathcal{P}_{\text{test}}$.
\\We are again interested in the behavior of the relative approximation error as a function of the number $N$ of basis functions used at the reduced order level.
We therefore let $N$ vary in a discrete set $\mathcal{N}$: for a fixed value of $N\in\mathcal{N}$, and for each $\theta_i$, $i=1, \dots, N_{\text{test}}$, we compute both the reduced solution $(\textbf{u}_N(t, \theta^i), p_N(t, \theta^i))$ and the corresponding high order solution $(\bm{u}_h(t, \theta^i), p_h(t, \theta^i))$. We calculate the $L^2$ relative error $\epsilon_{u,t_k}^{N, i}$ for the velocity and the relative error $\epsilon_{p,t_k}^{N, i}$ for the pressure at time $t_k$, by taking an average of these relatives error we obtain the mean approximation error $\epsilon_u^{N, i}$ for $\bm{u}$ and $\epsilon_p^{N, i}$ for $p$, for each $\theta_i\in\mathcal{P}_{\text{test}}$.
Finally we compute the average approximation errors $\overline{\epsilon}_u^N$ and $\overline{\epsilon}_p^N$ for every $N\in\mathcal{N}$, defined as:
\begin{equation*}
\overline{\epsilon}_u^N = \frac{1}{N_{\text{test}}}\sum_{i=1}^{N_{\text{test}}}\epsilon_u^{N, i}.
\end{equation*}

\begin{figure}
\centering
\includegraphics[scale=0.475]{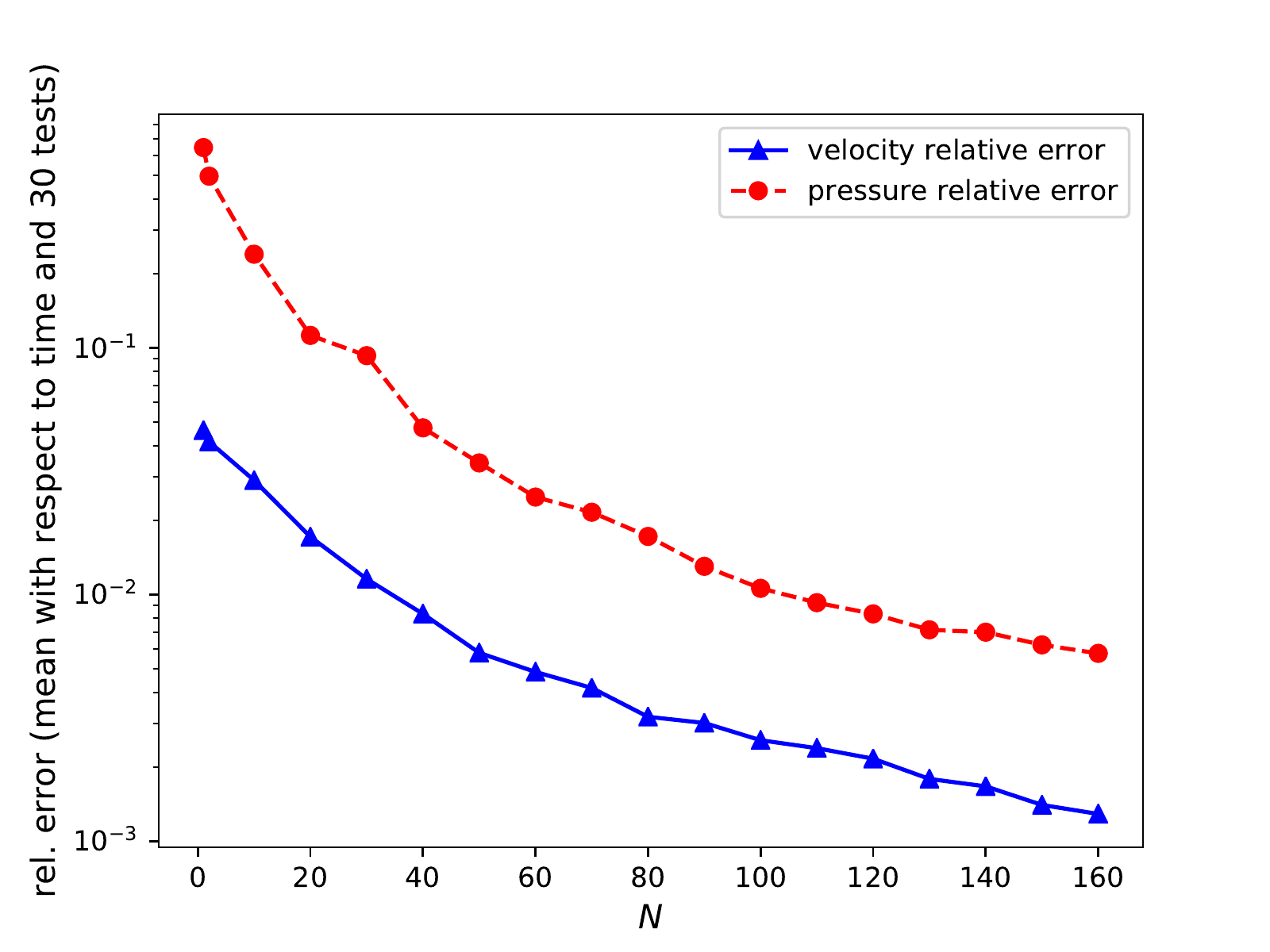}
\caption{Unsteady case: Mode dependent errors between high fidelity and reduced order approximation, with the supremizer enrichment.}
\label{unsteady: r_err_with_supremizer}
\end{figure}

Figure \ref{unsteady: r_err_with_supremizer} shows the relative approximation errors plotted against the number $N$ of basis functions used, with the supremizer enrichment at the reduced order level. 

\begin{figure}
\centering
\includegraphics[scale=0.25]{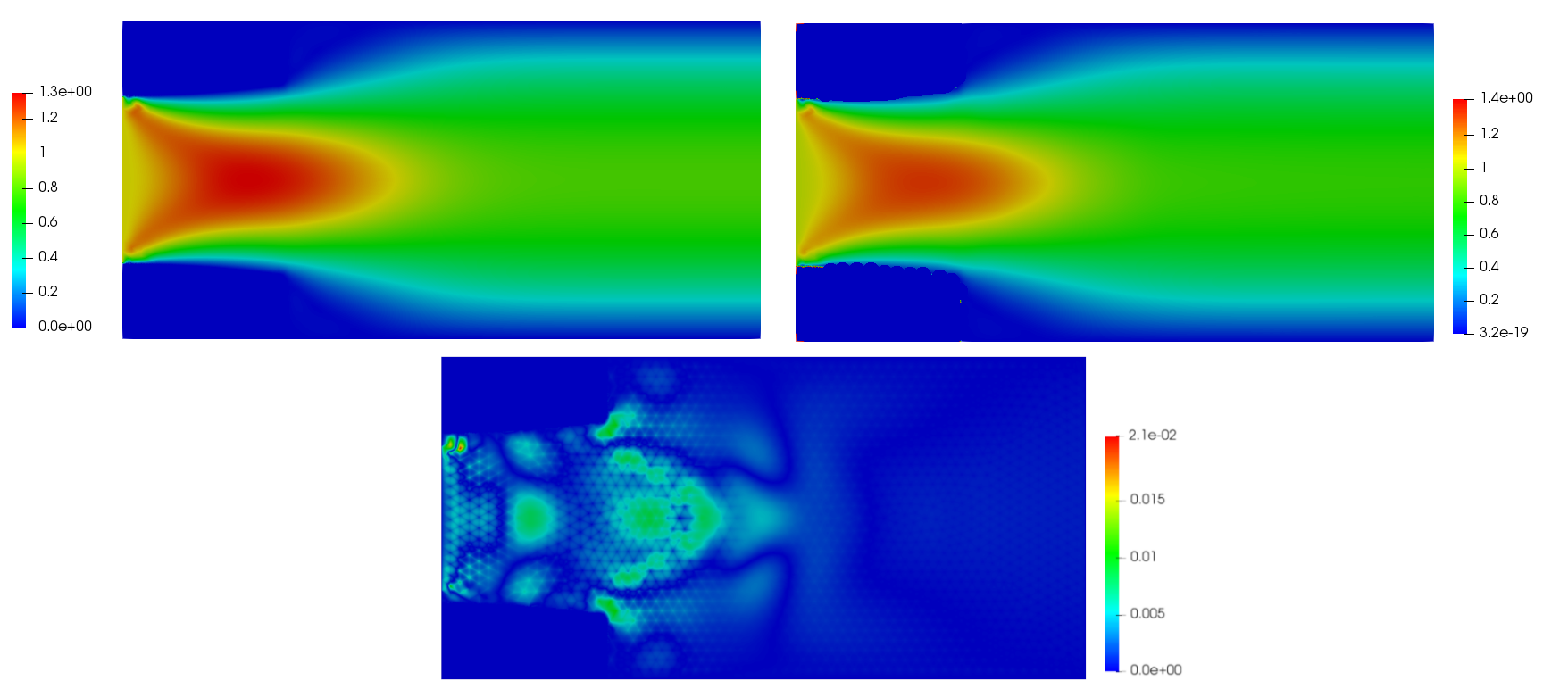}
\caption{Unsteady case: Cut geometry and high fidelity fluid velocity at final time $T=0.7$ for parameter $\theta=0.045406$ (left), reduced order solution for the same $\theta$ (right) and approximation error (bottom).}
\label{unsteady:fluid_velocity}
\end{figure}

\begin{figure}
\begin{center}
\includegraphics[scale=0.25]{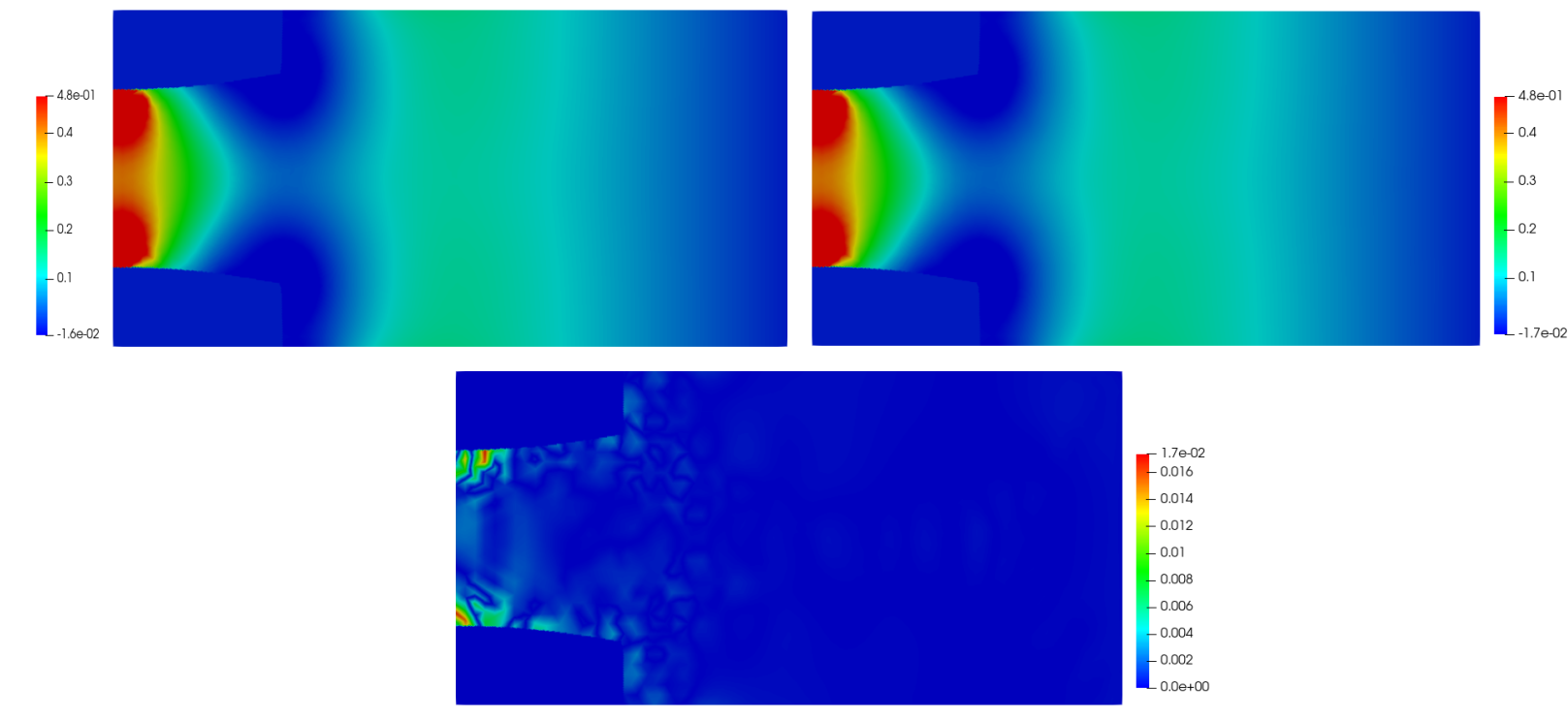}
\caption{Unsteady case: Cut geometry and the high fidelity pressure solution at final time at final time $T=0.7$ for parameter $\theta=0.050014$ (left), reduced order solution for the same $\theta$ (right) and approximation error (bottom).}
\label{unsteady:pressure_approx_error_with_s}
\end{center}
\end{figure}

Figure \ref{unsteady:pressure_approx_error_with_s} shows the approximation error for the pressure, for the same parameter value, with the supremizer enrichment. Figure \ref{unsteady:fluid_velocity} shows the approximation error for the fluid velocity $\bm{u}_f$ for a given value of the test parameter $\theta$, at the final time-step of the simulation.

\section{Unsteady Navier--Stokes: immersed obstacle}\label{immersed obstacle}
We now consider here a further test case of interest, namely the case of an obstacle immersed in a fluid.  For this test case we assume that $\mathcal{D}(\theta)$ represents a cylinder immersed in the fluid domain, and therefore we denote herein $\Gamma_{\theta} = \partial\mathcal{D}(\theta)$ the immersed boundary. 
\begin{figure}
\centering
\begin{tikzpicture}
\node[anchor=south west,inner sep=0] (image) at (0,0) {\includegraphics[scale=0.4]{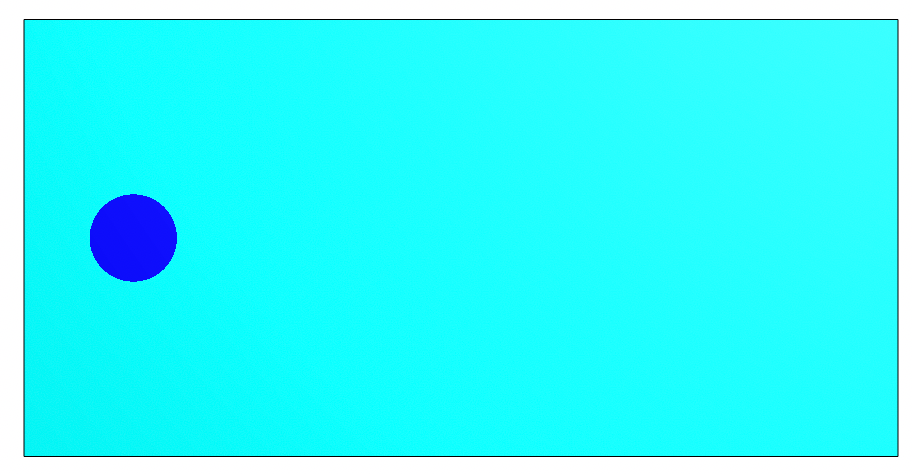}};
\begin{scope}[x={(image.south east)},y={(image.north west)}]
	\draw[black,thick] (0.07,0.3) node {$\Gamma_{in}$};
	\draw[black,thick] (0.88,0.3) node {$\Gamma_{out}$};
	\draw[black,thick] (0.5,0.9) node {$\Gamma_{D}$};
	\draw[black,thick] (0.5,0.1) node {$\Gamma_{D}$};
	\draw[black,thick] (0.23,0.5) node {$\Gamma_{\theta}$};
	\draw[white,thick] (0.15,0.5) node {\small{$\mathcal{D}(\theta)$}};
\end{scope}
\end{tikzpicture}
\caption{The physical domain of interest (light blue) and the levelset geometry (dark blue). The position of the immersed obstacle is determined by a geometrical parameter $\theta$. The immersed boundary is $\partial\mathcal{D}(\theta):=\Gamma_{\theta}$.}\label{cylinder domain}
\end{figure}
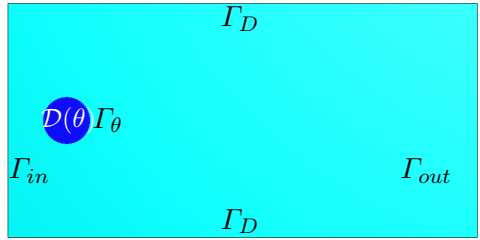

The physical domain over which the problem is formulated is depicted in Figure \ref{cylinder domain}.
\subsection{Strong formulation.}
The problem reads as follows: for every $t\in[0, T]$ and for every $\theta\in\mathcal{P}\subset\mathbb{R}$, find $\bm{u}(t;\theta)\colon\Omega(\theta)\mapsto\mathbb{R}^2$, $p(t; \theta)\colon\Omega(\theta)\mapsto\mathbb{R}$ such that:
\begin{equation}
\label{time dependent levelset Navier-Stokes}
\begin{cases}
\partial_t\bm{u}(\theta)-\mu\Delta \bm{u}(\theta) + \nabla p(\theta) + (\bm{u}(\theta)\cdot \nabla)\bm{u}(\theta) =  \bm{f}(\theta) \quad \text{in $\Omega(\theta)\times[0,T]$}, \\
\text{div}\bm{u}(\theta)=0 \quad \text{in $\Omega(\theta)\times[0,T]$},\\
\bm{u}(\bm{x},0;\theta)={\bm{u}}^0(\bm{x}, \theta) \quad \text{in $\Omega(\theta)$}.
\end{cases}
\end{equation}
The previous system is completed by the following boundary conditions: a prescribed inlet velocity $\bm{u}_{in}$ at the inlet boundary $\Gamma_{in}$, a no slip boundary condition on the immersed boundary $\Gamma_{\theta}$, a zero outflow condition on $\Gamma_{out}$ and a boundary condition $\bm{u}(\theta)\cdot\bm{n}_D=0$ on $\Gamma_D$.

\subsubsection{Geometrical parametrization}
The obstacle immersed in the fluid domain in our problem is a circle, defined through the \emph{time dependent} levelset function:
\begin{equation*}
\phi(x, y, \theta) = (x +1.5)^2 + (y -\theta)^2 - R^2,
\end{equation*}
where $\theta$ determines the position of the center of the cylinder in the domain, and $R$ is the radius of the circle.
\begin{figure}
\centering
\includegraphics[scale=0.75]{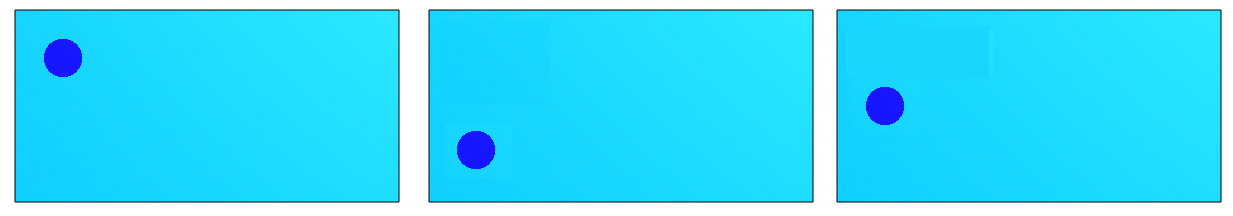}
\caption{Different levelset geometries: from left to right, $\theta=0.5$, $\theta=-0.5$, $\theta=0$.}\label{fig:lset_t_mesh}
\end{figure}
Figure \ref{fig:lset_t_mesh} shows the physical domain $\Omega(\theta)$ for different values of the parameter $\theta$: as we can see, changing the value of the geometrical parameter can produce a significant change in the physical domain of interest, and therefore this situation is particularly interesting for an unfitted discretization point of view, since adopting a standard discretization would require remeshing, or, alternatively, remapping the whole problem \eqref{time dependent levelset Navier-Stokes} onto a reference configuration, similarly to what is done for example with an Arbitrary Lagrangian Eulerian approach in fluid--structure interaction (see for example \cite{NoBaRo19, BallarinRozza2016, Richter, Wang2018}).

\subsubsection{Weak formulation and time discretization}\label{immersed obstacle weak formulation}
We now want to state the weak formulation of the original problem after discretization in space and after having applied a time stepping scheme. As far as the time discretization concerns, we employ the time stepping scheme adopted in Section \ref{Unsteady NS}: we discretize the time interval $[0, T]$ in sub-intervals $(t^n, t^{n+1}]$ of measure $\tau_{n+1}=t^{n+1}-t^n$, for $n=0,\hdots, N_t-1$. Also in this case, we decide to treat the boundary conditions in the following way: the boundary conditions on $\Gamma_D$ and on $\Gamma_{out}$, as well as the condition on the inlet profile, are imposed strongly; only the boundary condition on the immersed boundary $\Gamma_{\theta}$ is imposed weakly, thanks to the Nitsche method.
For the space discretization we therefore use the discrete spaces introduced in Section \ref{Time discretization}, and which we recall briefly:
\begin{equation*}
\begin{split}
V_{h,2}(\theta) &:= \{\bm{v}_h \in (C_0(\Omega^{*}_h(\theta)))^2 \colon \bm{v}_h |_{T} \in (\mathcal{P}^2(T))^2, \quad\forall T \in \mathcal{I}_h(\theta)\},\\
V_{h,2}^D(\theta)&:=\{\bm{v}_h\in V_{h,2}(\theta) \text{ such that } \bm{v}_h=\bm{u}_{in}\text{ on }\Gamma_{in} \text{ and }\bm{v}_h\cdot\bm{n}_D=0\text{ on }\Gamma_D\},\\
V_{h,2}^0(\theta)&:= \{\bm{v}_h\in V_{h,2}(\theta) \text{ such that } \bm{v}_h=\bm{0}\text{ on }\Gamma_{in} \text{ and } \bm{v}_h\cdot\bm{n}_D=0\text{ on }\Gamma_D\},\\
Q_{h,1}(\theta) &:= \{q_h \in C_0(\Omega^{*}_h(\theta))\colon q_h|_{T} \in \mathcal{P}^1(T), \quad\forall T \in \mathcal{I}_h(\theta)\}.
\end{split}
\end{equation*}
For the weak formulation we restore to the one presented in Section \ref{Time discretization}.
\subsection{Proper Orthogonal Decomposition and online system}\label{immersed obstacle POD}
Also in this case, once we have computed the snapshots $(\bm{u}_h(t^i;\theta), p_h(t^i;\theta))$, for all $\theta$ in a discrete parameter training set $\mathcal{P}_{train}$ and for $i=0,\dots, N_t$, we rely on a natural smooth extension of the snapshots, and on a lifting of the inlet boundary condition, in order to obtain solutions that are defined on the whole background mesh $\hat{\mathcal{I}}_h$. Again, we use a supremizer enrichment technique in order to obtain stable approximations of the pressure in the online step, and we use a lifting function to treat the non--homogeneous Dirichlet condition at the inlet boundary $\Gamma_{in}$. Also here, in order to make the notation more light, we will omit the lifting function, with the understanding that it has been used for implementation purposes.
After having created the \emph{parameter independent} reduced basis functions $\{\bm{\Phi}_i^{u,s}\}_{i=1}^{2N}$ and $\{\Phi_i^p\}_{i=1}^N$ for the fluid velocity and fluid pressure respectively, the online algebraic system that we are going to solve, for every $\theta\in\mathcal{P}$, is the one presented in Equation \eqref{eq:system_linear_reduced}:
\begin{equation*}
\begin{bmatrix}
{\bm{L}_{u,s}^T\hat{\bm{M}}\bm{L}_{u, s}} & \bm{0}\\
\bm{0} & \bm{0}
\end{bmatrix}
U_N^{n+1}(\theta)
+
\tau_{n+1}
\hat{R}_N(U_N^{n+1}(\theta); \theta)
=
\begin{bmatrix}
{{\bm{L}_{u,s}^T\hat{\bm{M}}\bm{L}_{u, s}}} & \bm{0}\\
\bm{0} & \bm{0}
\end{bmatrix}
U_N^n(\theta).
\end{equation*}

\subsection{Numerical results}\label{immersed obstacle numerical results}
We now present some numerical results for the test case of the immersed obstacle. The cylinder has a radius $R=0.2$ $cm$, and the center C of the cylinder has coordinates $(x_c, y_c)=(-1.5, \theta)$, with $\theta\in[-0.65, 0.65]$. The background domain $\mathcal{R}$ is a rectangle of coordinates $(-2, -1)$ and $(2, 1)$. The fluid viscosity is $\mu=0.05$, the fluid density is $
1$. The mesh size is $h_{max}=0.07$, and the timestep used for the discretization is $\tau=\frac{h_{max}}{6}$; the final time of the simulation is $T=0.7$ $s$. The inlet velocity profile is $u_{in}=(1, 0)$ $m/s$. In this case, $N_{train}=200$ is the number of parameters $\theta^i$ that we take to compute the snapshots for the POD, whereas $N_{test}=30$ is the number of parameters $\theta^i$ for which we compute the online solution.
\begin{figure}
\centering
\includegraphics[scale=0.15]{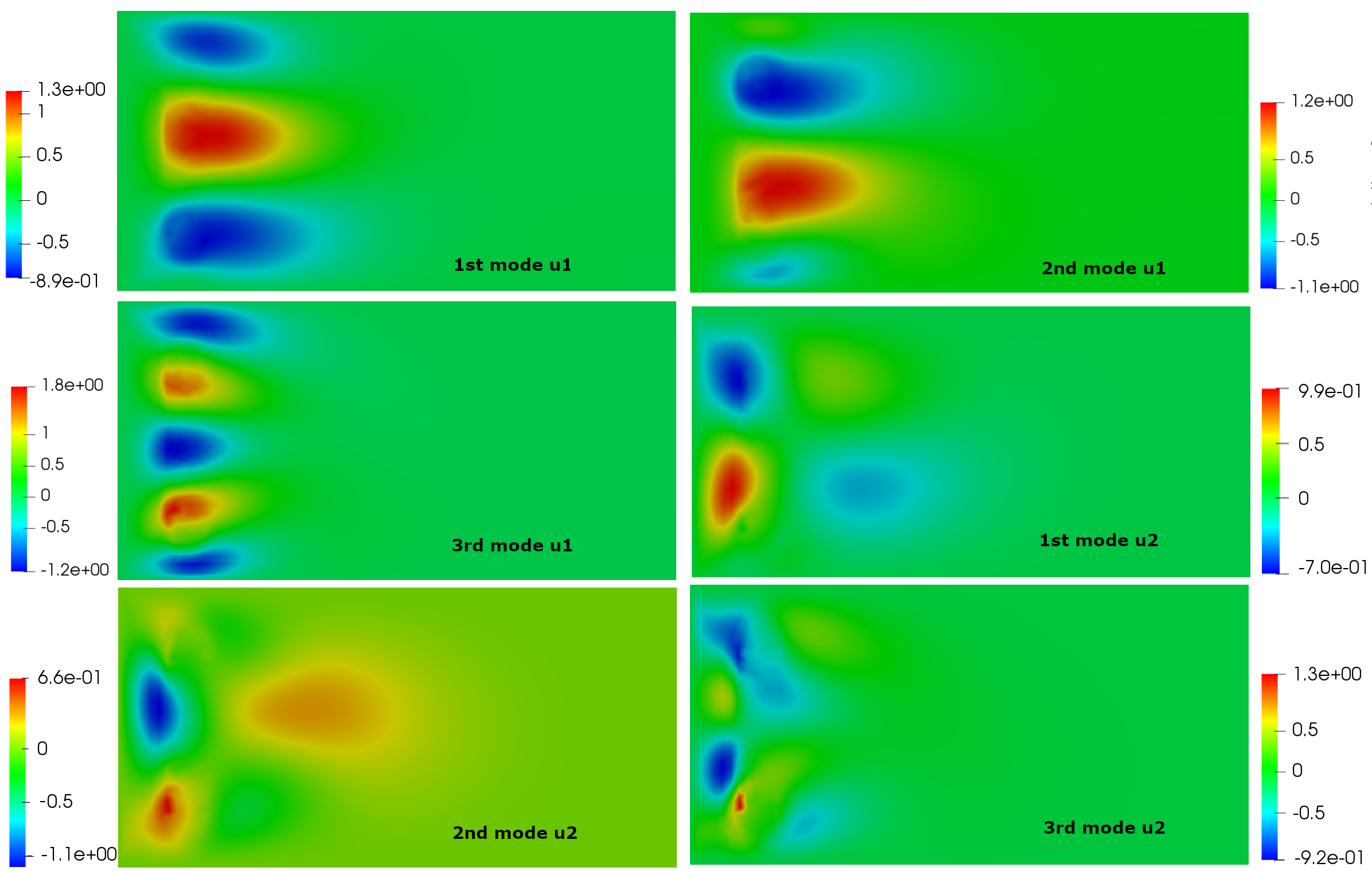}
\caption{First three modes for the $x$--component of the velocity and for the $y$--component. }\label{velocity_modes}
\end{figure}
\begin{figure}
\centering
\includegraphics[scale=0.15]{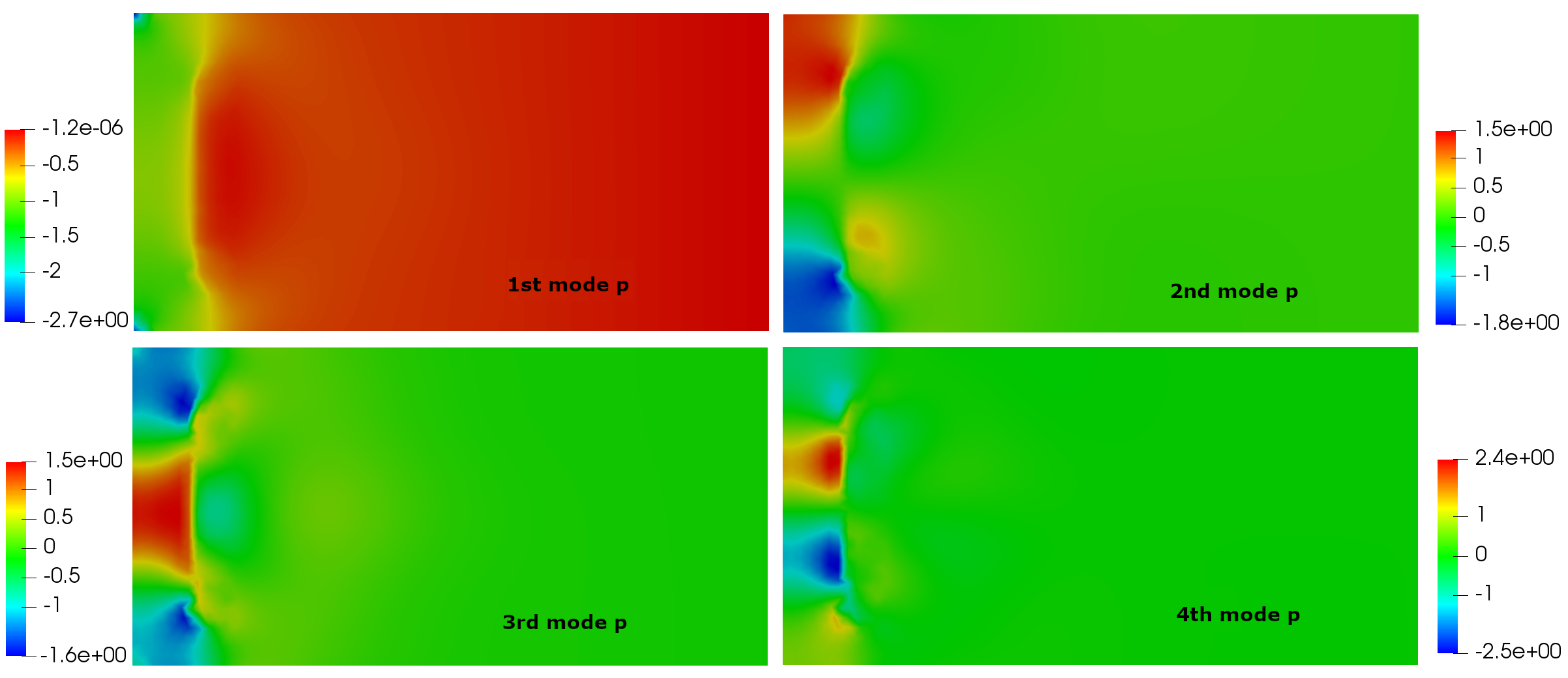}
\caption{First four pressure modes}\label{pressure_modes}
\end{figure} 
%

\begin{figure}
\centering
\includegraphics[scale=0.2]{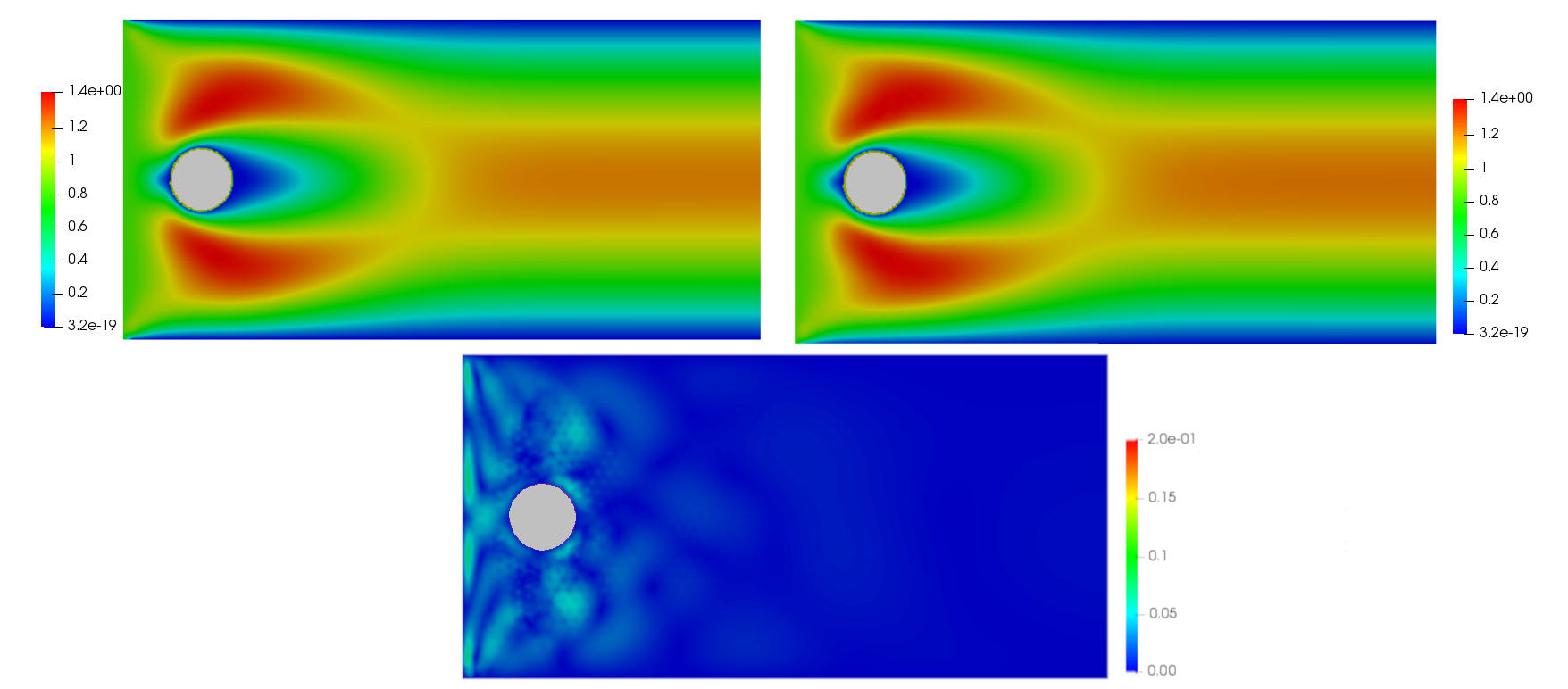}
\caption{Fluid velocity approximation (top left), reduced order approximation (top right) and approximation error (bottom), for $\theta=0$, at time $t=0.7$ s.}\label{velocity0}
\end{figure}
\begin{figure}
\centering
\includegraphics[scale=0.2]{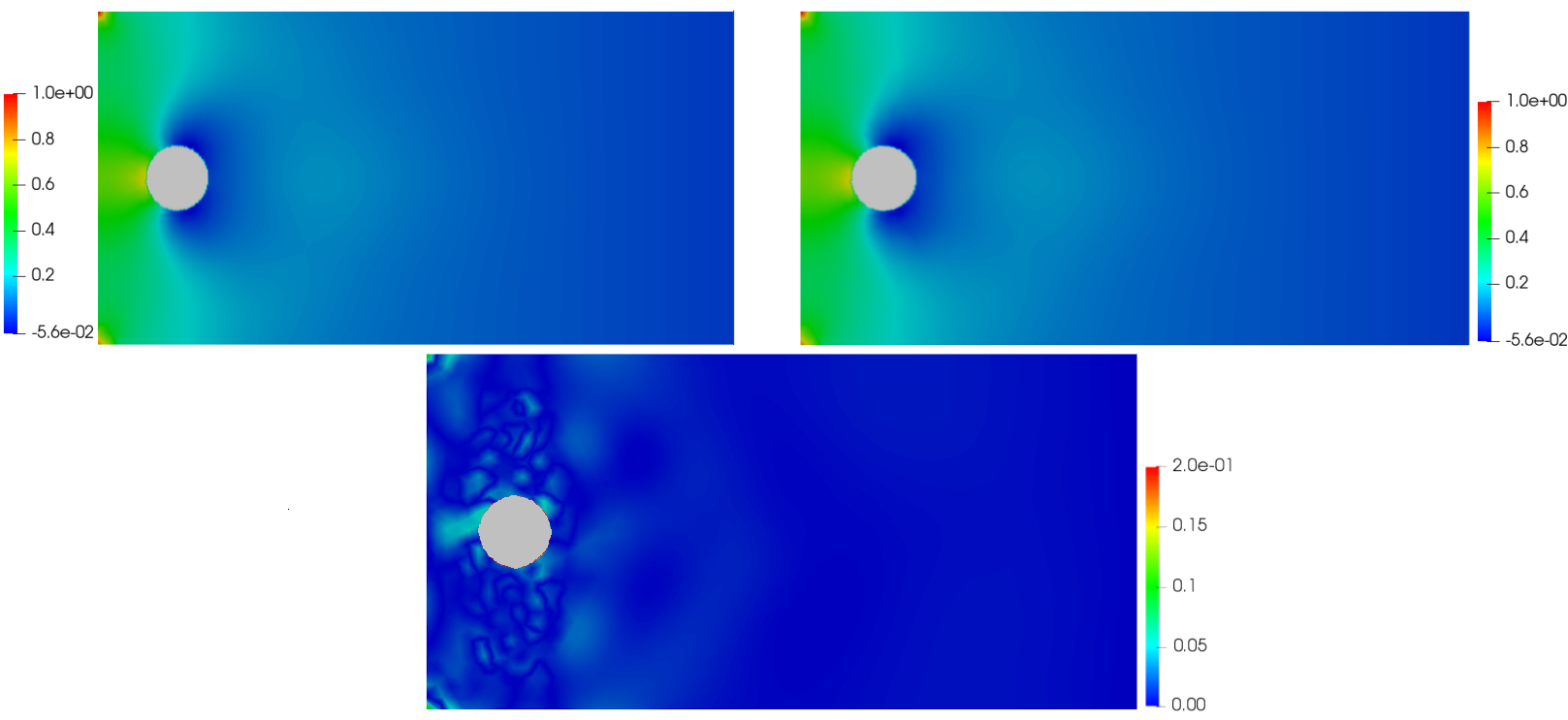}
\caption{Pressure approximation (top left), reduced order approximation (top right) and approximation error (bottom), for $\theta=0$, at time $t=0.7$ s.}\label{pressure0}
\end{figure} 
\begin{figure}
\centering
\includegraphics[scale=0.2]{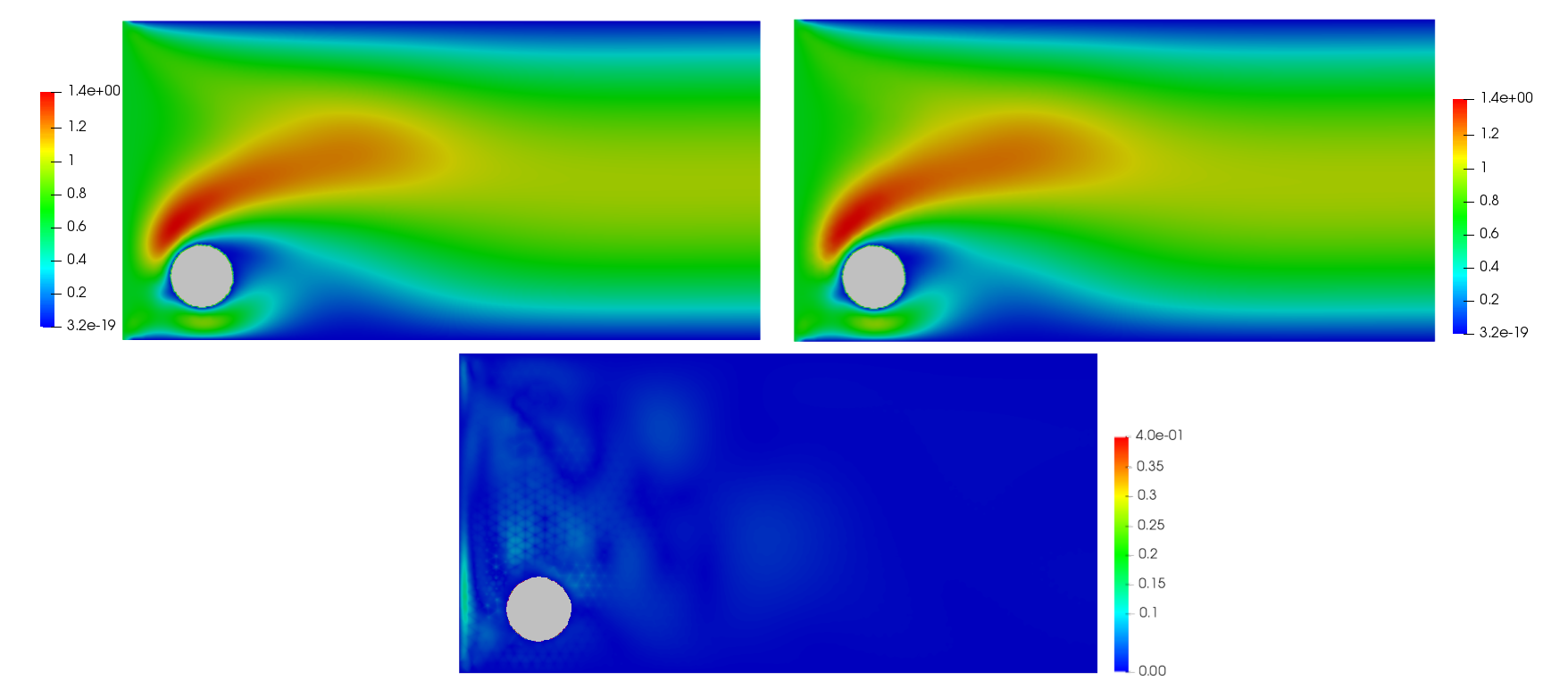}
\caption{Fluid velocity approximation (top left), reduced order approximation (top right) and approximation error (bottom), for $\theta=0.6$, at time $t=0.7$ s.}\label{velocity1}
\end{figure}
\begin{figure}
\centering
\includegraphics[scale=0.2]{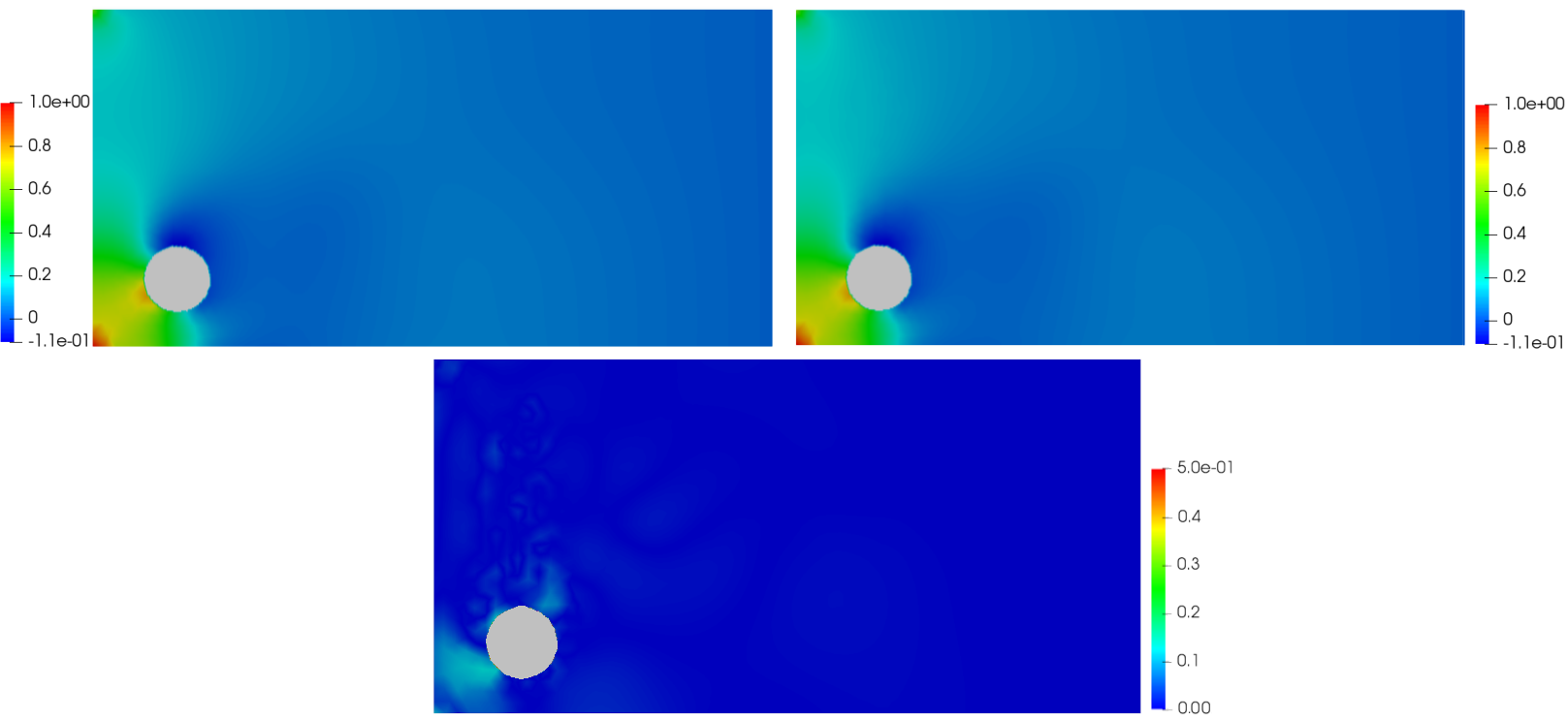}
\caption{Pressure approximation (top left), reduced order approximation (top right) and approximation error (bottom), for $\theta=0.6$, at time $t=0.7$ s.}\label{pressure1}
\end{figure} 
\begin{figure}
\centering
\includegraphics[scale=0.2]{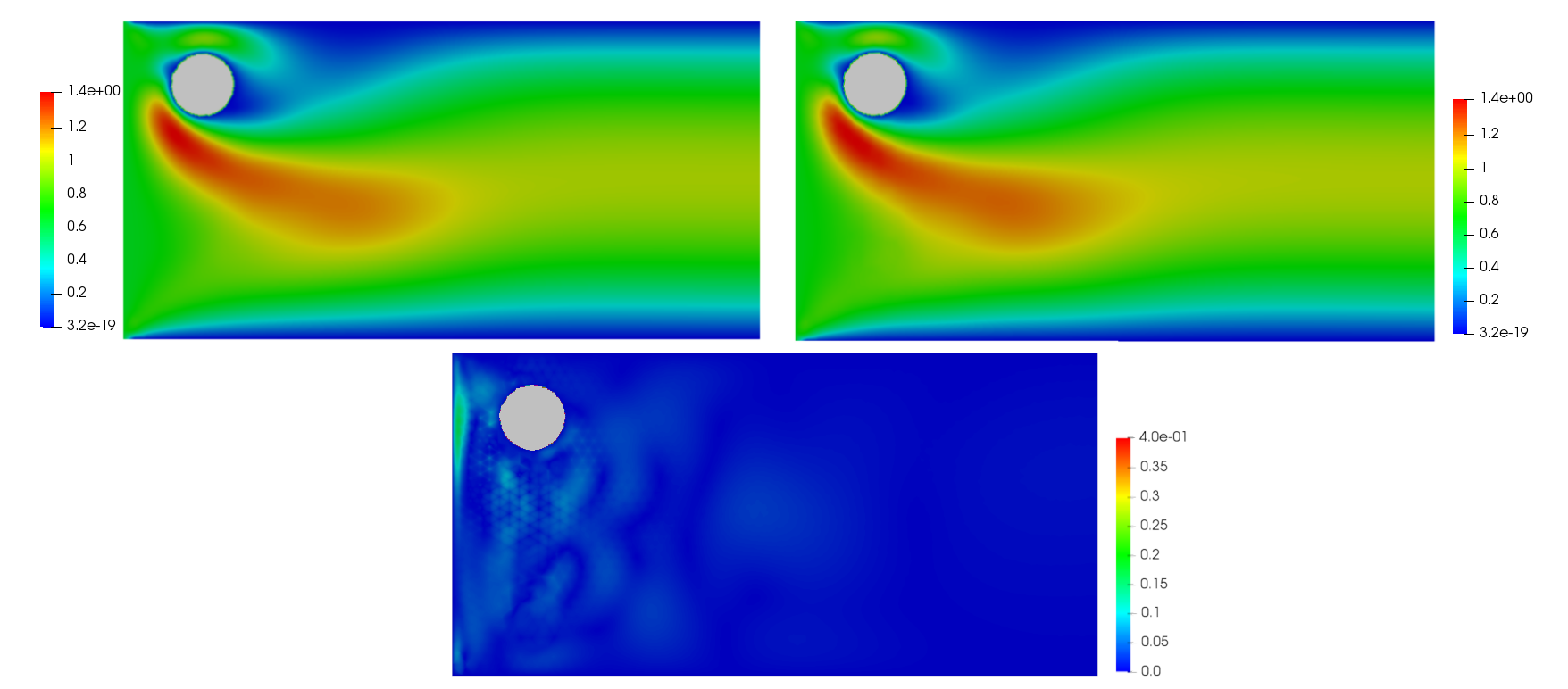}
\caption{Fluid velocity approximation (top left), reduced order approximation (top right) and approximation error (bottom), for $\theta=-0.6$, at time $t=0.7$ s.}\label{velocity2}
\end{figure}
\begin{figure}
\centering
\includegraphics[scale=0.2]{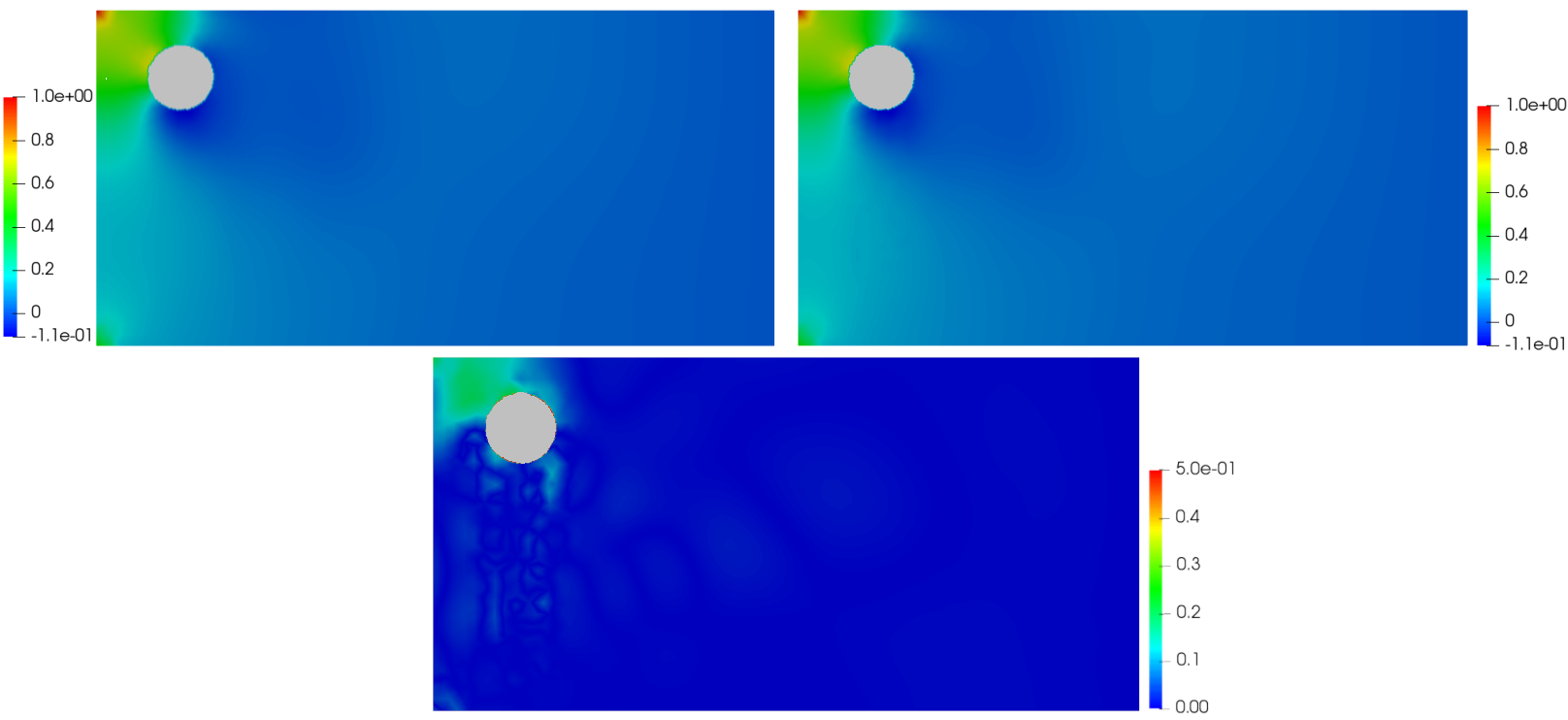}
\caption{Pressure approximation (top left), reduced order approximation (top right) and approximation error (bottom), for $\theta=-0.6$, at time $t=0.7$ s.}\label{pressure2}
\end{figure} 

\begin{figure}
\centering
\includegraphics[scale=0.5]{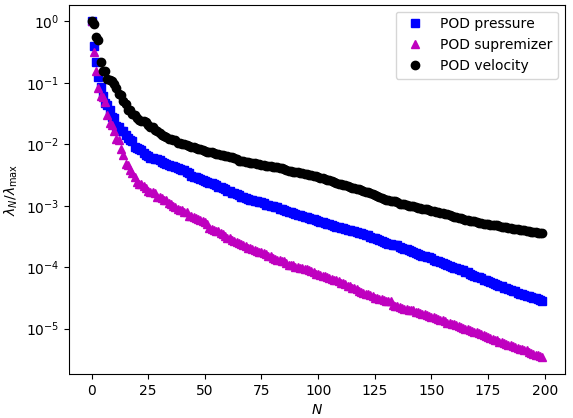}
\caption{Immersed obstacle: decay of the first $200$ eigenvalues returned by the POD on the fluid velocity (black), the pressure (magenta) and the supremizer (blue).}\label{eigs cylinder}
\end{figure}

\begin{figure}
\centering
\includegraphics[scale=0.5]{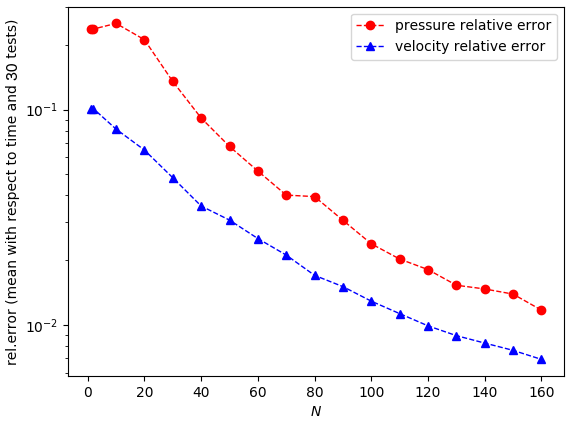}
\caption{Immersed obstacle: mean error behaviour, according to different number of basis functions used in the online phase, for the fluid velocity (blue) and the pressure (red). The results were obtained using the supremizer enrichment technique.}\label{mean error cylinder}
\end{figure}

Figures \ref{velocity0}--\ref{pressure2} represent the fluid velocity and the fluid pressure, at the last time-step of the simulation, for three different values of the geometrical parameter $\theta$. As we can see, the change in the position of the immersed obstacle, according to the parameter $\theta$, is significant; nonetheless, with an unfitted discretization we were able to obtain very good results, in terms of both velocity and pressure approximation. The use of the supremizer enrichment technique helps to obtain a more accurate results for the reduced order fluid pressure. In Figure \ref{eigs cylinder} we can see the behaviour of the first eigenvalues returned by a POD on the fluid velocity, the fluid pressure and the fluid supremizer; Figure \ref{mean error cylinder} shows the behaviour of the mean approximation error, plotted against the number $N$ of basis functions used for the online approximation, and for $30$ values of the parameter $\theta$ in the test sample.
We remark that these results have been obtained without the implementation of any snapshot transportation technique (see for example \cite{KaBaRO18, NoBaRo19, CaMaSta19}): the analysis of the influence of this additional feature on the overall quality of the approximation at the reduced order level is left for a future work.

\section{Conclusions}\label{sec:conclusions}
\noindent In this work we have introduced a POD--Galerkin ROM approach for geometrically parametrized two dimensional Navier--Stokes equations, both in the steady and in the unsteady case. 
The procedure that we have proposed shows many of the advantages that characterize CutFEM and reduced order methods. {First of all, by choosing an unfitted mesh approach, we have shown that it is possible to work with geometries that can potentially change significantly the shape of (part of) the domain, as we can see from the examples in Figure \ref{levelset_image} and Figure \ref{fig:lset_t_mesh}. By employing an unfitted CutFEM approach at the full order level, we can let the geometrical parameter $\theta$ vary in a large interval of values: this helps to overcome one of the limitations of the standard Finite Element discretization, where a re--meshing would have been needed.} 

At the reduced order level, {an unfitted mesh discretization has been combined with a snapshot extension, in order to be able to work with finite dimensional spaces that are parameter--independent: the importance of this aspect is related to the choice of performing a Proper Orthogonal Decomposition. In addition, the implementation of a supremizer enrichment technique allows us to obtain a stable approximation of the fluid pressure. Combined together, the unfitted CutFEM discretization and the Reduced Basis Method represent a very powerful tool, which here allowed us to investigate time--dependent, nonlinear fluid flows problems.}
To conclude, with this work we started to prepare the basis for a CutFEM-RB procedure that, thanks to a Galerkin projection and with the use of a supremizer enrichment, will ideally be used to obtain accurate approximations of solutions of very complex problems, such as fully coupled multiphysics problem, where large displacements of the structure may occur.
Therefore, as future perspectives, we first would like to extend the procedure presented here to nonlinear, time--dependent problems, where the levelset geometry is time--dependent (i.e. the motion of the levelset geometry is time--dependent, and this motion is known a priori); then, we would like to move to a more general case, where on the contrary the motion of the levelset geometry is an unknown of the problem. As a next step and as final goal therefore, we would like to test the performance of this approach with time dependent Fluid--Structure Interaction problems, geometrically shallow water flows as well as with phase flow Navier-Stokes systems. Furthermore, from the model reduction point of view, we will pursue further developments in hyper-reduction techniques \cite{Xiao20141,BARRAULT2004667,Carlberg2013623, stabile_geo_} tailored for unfitted discretizations.

\section*{Acknowledgements}
We acknowledge the support by European Union Funding for Research and Innovation -- Horizon 2020 Program -- in the framework of European Research Council Executive Agency: Consolidator Grant H2020 ERC CoG 2015 AROMA-CFD project 681447 "Advanced Reduced Order Methods with Applications in Computational Fluid Dynamics'' (PI Prof. Gianluigi Rozza).
We also acknowledge the INDAM-GNCS project "Tecniche Numeriche Avanzate per Applicazioni Industriali''.
The first author has received funding from the Hellenic Foundation for Research and Innovation (HFRI) and  the  General  Secretariat  for  Research  and  Technology (GSRT), under  grant agreement No[1115], the "First Call for H.F.R.I. Research Projects to support Faculty members and Researchers and the procurement of high-cost research equipment" grant 3270 and the National Infrastructures for Research and Technology S.A. (GRNET S.A.) in the National HPC facility - ARIS - under project ID pa190902. The second author acknowledges also the support of the Austrian Science Fund (FWF) project F65 "Taming complexity in Partial Differential Systems" and the Austrian Science Fund (FWF) project P 33477.\\
Numerical simulations have been obtained with the extension \emph{ngsxfem} of \emph{ngsolve} \cite{ngsxfem,ngsolve,schoeberl} for the high fidelity part, and RBniCS \cite{rbnics} for the reduced order part. We acknowledge developers and contributors of each of the aforementioned libraries.
 %
\bibliographystyle{amsplain}
\bibliography{literature}
\end{document}